\documentclass{amsart} \usepackage{amssymb,amsfonts,amscd}
\newtheorem{theorem}{Theorem}[section]
\newtheorem{lemma}[theorem]{Lemma}
\newtheorem{corollary}[theorem]{Corollary}
\newtheorem{proposition}[theorem]{Proposition}
\theoremstyle{remark}

\theoremstyle{definition}
\newtheorem{definition}[theorem]{Definition}

\numberwithin{equation}{section} \makeatother

\DeclareMathOperator{\Kdb}{{\mathbb K}}

\DeclareMathOperator{\Cdb}{{\mathbb C}}
\DeclareMathOperator{\Rdb}{{\mathbb R}}
\DeclareMathOperator{\Ddb}{{\mathbb D}}

\DeclareMathOperator{\Ndb}{{\mathbb N}}

\DeclareMathOperator{\cA}{{\mathcal A}}
\def\alp{{\alpha}}
\def\bN{{\ensuremath{\mathbb N} }}
\def\bC{{\ensuremath{\mathbb C} }}

\def\bR{{\ensuremath{\mathbb R} }}
\def\nm#1{{\left\Vert #1 \right\Vert}}
\def\Nm{{\nm\cdot}}
\def\veps{{\varepsilon}}
\newcommand\bdfn{\begin{defn}}
\newcommand\blem{\begin{lemma}}
\newcommand\elem{\end{lemma}}
\newcommand\bcor{\begin{corollary}}
\newcommand\ecor{\end{corollary}}
\newcommand\bthm{\begin{theorem}}
\newcommand\ethm{\end{theorem}}
\newcommand\edfn{\end{defn}}
\def\gam{{\gamma}}

\def\loc{{\textrm {loc}}}
\def\del{{\delta}}
\newif\ifrough
\roughfalse
\def\reff#1{\ifrough{\ref {#1}=#1}\else{\ref {#1}}\fi}
\def\l#1!{{\ifrough{\textrm{=#1\ }\label{#1}}\else{\label{#1}}\fi}}

\def\summ#1#2#3{{\sum_{#1=#2}^{#3}}}
\def\reflem#1{{Lemma \reff{#1}}}

\def\refthm#1{{Theorem \reff{#1}}}
\def\refcor#1{{Corollary \reff{#1}}}
\def\ni{{\noindent}}
\begin{document}

\title[Operator algebras with cai]{Operator algebras with contractive approximate identities} \author[D. P. Blecher]{David P. Blecher}
\author[C. J. Read]{Charles John Read}

\address{Department of Mathematics, University of Houston, Houston, TX
77204-3008}
 \email[David P.
Blecher]{dblecher@math.uh.edu}
\address{Department of Pure Mathematics,
University of Leeds,
Leeds LS2 9JT,
England}
 \email[Charles John Read]{read@maths.leeds.ac.uk}
\date{11/13/2010.  Revision of 2/18/2011.}
\thanks{*Blecher was partially supported by a grant  from
the National Science Foundation.}
\begin{abstract}
We give several applications of a recent theorem  of the second author,
which solved a conjecture  of the
first author with Hay and
Neal, concerning  contractive approximate identities; and another of
Hay from the theory of noncommutative peak sets, thereby putting the
latter theory on a much firmer foundation. From this theorem it
emerges there is a surprising amount of positivity present in any
operator algebras with contractive approximate identity. We exploit
this to generalize several results previously available only for
$C^*$-algebras, and we give many other applications.
\end{abstract}

\maketitle

\section{Introduction}

 An {\em operator algebra} is a closed subalgebra of $B(H)$, for a
Hilbert space $H$. We recall that by a theorem due to  Ralf Meyer,
every operator algebra $A$ has a unique unitization $A^1$ (see
\cite{Mey} or \cite[Section 2.1]{BLM}). Below $1$ always refers to
the identity of $A^1$ if $A$ has no identity.   We are mostly
interested in operator algebras with  contractive approximate
identities (cai's).  We also call these {\em approximately unital}
operator algebras. In our paper we give several applications of the
following recent result, which was prompted by, and solves, a
question on p.\ 351 of \cite{BHN}:

\begin{theorem} \label{read} {\rm \cite{Read} \ }  An operator algebra
with a cai, has a cai $(e_t)$ with $\Vert 1 - e_t \Vert \leq 1$, and even
with $\Vert 1 - 2 e_t \Vert \leq 1$, for all $t$.
 \end{theorem}

This result draws attention to the set of operators $x$ in an
operator algebra $A$ satisfying $\Vert 1 - x  \Vert \leq 1$.  We
denote this set by ${\mathfrak F}_A$; it will play a role for us
very much akin to the role of the positive cone in a $C^*$-algebra.
This surprising claim will be justified at many points in our paper,
but the reader can begin to see this by considering the following
fact: a linear map $T : A \to B$ between $C^*$-algebras or operator
systems is completely positive in the usual sense iff there is a
constant $C
> 0$ such that $T({\mathfrak F}_{A}) \subset C {\mathfrak F}_{B}$,
and similarly at the matrix levels (see Section 8). Indeed we use
Theorem \ref{read} to see that there is a remarkable amount of
positivity present in any operator algebra with cai.  We exploit
this, and various properties of operators in ${\mathfrak F}_A$, to
generalize several results previously available only for
$C^*$-algebras.  Many of the applications which we give are to the
structure theory of operator algebras.   Some of these advances are
mentioned in more detail in the next paragraph.
 We recall that a classical principle is to study a
ring or algebra $A$ in terms of its ideals, both two-sided and one-sided.  Unfortunately, not much is
known about general closed ideals in $A$, even for common examples
of function algebras, and so we focus on the {\em r-ideals}
(right ideals with a left cai) and {\em $\ell$-ideals} (left ideals with a
right cai).
As proved in \cite{BHN}, these objects are in an inclusion-preserving,
bijective correspondence with each other, and also with the
{\em hereditary subalgebras} (or {\em HSA's}; defined below).
HSA's are frequently more useful, for example in $C^*$-algebra theory,
because they are
more symmetrical objects, and because many important properties pass to
HSA's \cite{Bla}.

The layout of our paper is as follows:  At the end of  Section 1 we give some
quick consequences of Theorem \ref{read}.
The long Section 2 contains a number of facts about
${\mathfrak F}_A$, and uses these together with Theorem \ref{read}
to give many applications to the structure  of
operator algebras.  For example, one theme of our paper is how
cai's may be built.  We dissolve the remaining
mysteries concerning r-ideals by showing  how they all arise.
The separable r-ideals in a operator algebra are precisely the
subspaces $\overline{xA}$, for an element $x \in {\mathfrak F}_A$ which we may
select to be as close as we like to a positive norm 1 operator. The
nonseparable r-ideals are limits of increasing nets of such subspaces
$\overline{xA}$.   Similarly for the matching class of $\ell$-ideals, or HSA's.
 Other sample results: we
show that as in the $C^*$-algebra case, any separable operator
algebra with cai has a countable cai consisting of mutually
commuting elements; and we prove
a noncommutative Urysohn lemma.  In Section 3 we study
the pseudo-invertible (sometimes called `generalized invertible')
elements in operator algebras, a topic connected to the notion
of `well supported' elements.
This topic is also very intimately  connected
to the
question of when a `principal ideal' $xA$ is already closed.
In Section 4 we study
operator algebras possessing no r-ideals or HSA's.  We also give several
interesting examples of such algebras.  In Section 5 we display a radical,
approximately unital operator algebra which is an integral domain, and whose
ideal structure can be completely determined.   Hence
this is an excellent example against which to
test certain
conjectures concerning the structure theory of operator algebras.
 In Section 6 we consider preimages of r-ideals, HSA's, etc.
In Section 7 we describe some other interesting constructions of
r-ideals in operator algebras.  In the final Section 8 we introduce
and study a notion of {\em completely positive maps} between general
operator algebras, or between unital operator spaces, and give an
Arveson type extension theorem, and a Stinespring type
characterization,  for such maps.

We remark that most of our results apply immediately to function algebras,
that is to uniformly closed subalgebras of $C(K)$ spaces, since these
are special cases of operator algebras.  We will not take the time
to point these out, although some of these applications are new.

We now state our notation, and some facts. We refer the reader to
\cite{BLM} for additional background on operator algebras, and for
some of the details and notation below.   For us a {\em projection}
is always an orthogonal projection, and an {\em idempotent} merely
satisfies $x^2 = x$. If $X, Y$ are sets, then $XY$ denotes the
closure of the span of products of the form $xy$ for $x \in X, y \in
Y$.   We write $X_+$ for the positive operators that happen to
belong to $X$. Returning to the unitization, if $A$ is a nonunital
operator algebra represented (completely) isometrically on a Hilbert
space $H$ then one may identify $A^1$ with $A + \Cdb I_H$.   The
second dual $A^{**}$ is also an operator algebra with its (unique)
Arens product, this is also the product inherited from the von Neumann
algebra $B^{**}$ if
$A$ is a subalgebra of a $C^*$-algebra $B$.  Meets and joins in
$B^{**}$ of projections in $A^{**}$ remain in $A^{**}$, as can be
readily seen for example by inspecting some of the classical
formulae for meets and joins of Hilbert space projections,
or by noting that these meets and joins occur in the biggest
von Neumann algebra contained inside $A^{**}$. Note that
$A$ has a cai iff $A^{**}$ has an identity $1_{A^{**}}$ of norm $1$,
and then $A^1$ is sometimes identified with $A + \Cdb 1_{A^{**}}$.
In this case the multiplier algebra $M(A)$ is identified with the
idealizer of $A$ in $A^{**}$ (that is, the set of elements
$\alpha\in A^{**}$ such that $\alpha A\subset A$ and $A
\alpha\subset A$).
 It can also be viewed as the
idealizer of $A$ in $B(H)$, if the above representation  on $H$ is  nondegenerate.
If  $A$ is unital then $M(A) = A$, and in this case
we often assume that $1_A = I_H$.

Let $A$ be an operator algebra.  The set
${\mathfrak F}_A = \{ x \in A : \Vert 1 - x \Vert \leq 1 \}$ equals
$\{ x \in A : \Vert 1 - x \Vert = 1 \}$ if $A$ is nonunital, whereas
if $A$ is unital then ${\mathfrak F}_A = 1 + {\rm Ball}(A)$.
If $x \in {\mathfrak F}_A$ then the numerical range of $x$ is contained
in the closed disk of center $1$ and radius $1$, and in particular
is in the right half plane (that is, $x$ is {\em accretive}).
Clearly  $x$ is a sectorial
 operator.  See \cite{sect} for more information on
sectorial and
accretive operators and their functional calculus.
Note that $x \in {\mathfrak F}_A$ iff $x x^* \leq x + x^* = 2 \, {\rm Re}(x)$,
and iff $x^* x  \leq x + x^* = 2 \, {\rm Re}(x)$.
If $A$ is a closed subalgebra of an operator algebra $B$
then it is easy to see, using the
uniqueness of the unitization, that ${\mathfrak F}_A = A \cap {\mathfrak F}_B$. We write $\frac{1}{2} {\mathfrak F}_A$ for $\{ x \in A : \Vert 1 - 2x \Vert \leq 1 \}$.
  We remark that the condition $||1- 2x|| \leq 1$ implies both $||x|| \leq 1$ and
$||1-x|| \leq 1$.  In much of our paper, where we have ${\mathfrak
F}_A$ it probably would be preferable to employ $\frac{1}{2}
{\mathfrak F}_A$ instead. However since in these occurrences it will
not matter technically, we use the simpler notation.

We recall that an {\em r-ideal} is a right ideal with a left cai, and an {\em $\ell$-ideal} is a left ideal with a right cai.
We say that an operator algebra $D$ with cai, which is a subalgebra of
another operator algebra $A$, is an HSA (hereditary subalgebra)
of $A$, if $DAD \subset D$.
For the theory of HSA's see \cite{BHN}.  These objects are in an order preserving,
bijective correspondence with the r-ideals in $A$, and also with the
{\em open projections} $p \in A^{**}$, by which we mean that there
is a net $x_t \in A$ with $x_t = p x_t p \to p$ weak*.  These are
also the open projections $p$ in the sense of Akemann \cite{Ake,Ake2} in $B^{**}$, where $B$ is a $C^*$-algebra containing $A$, such that
$p \in A^{\perp \perp}$.   The complement (`perp') of an open projection is called a {\em closed projection}.  We spell out some of the correspondences above:
if $D$ is a HSA in $A$, then $J = DA$ is the matching  r-ideal.
The weak* limit of a cai for $D$, or of a left cai for $J$, is
an open projection, and is called the {\em support projection} of $J$ or $D$.  Conversely, if $p$ is an open projection in $A^{**}$, then
$pA^{**} \cap A$ and $pA^{**}p \cap A$ is the matching r-ideal and HSA pair
in $A$.  We also mention that suprema (resp.\ infima) of open
(resp.\ closed) projections  in $A^{**}$, remain in $A^{**}$,
by the fact mentioned two paragraphs earlier about meets and joins,
together with the $C^*$-algebraic case of these facts \cite{Ake,Ake2}.

The {\em peak} and {\em $p$-projections} studied in
 \cite{Hayth} and \cite{BHN}, are certain closed projections which generalize the important notions
of {\em peak sets} and {\em $p$-sets} from the theory of function spaces.
We recall that a peak set for a unital
space $A$ of continuous functions
on a compact set $K$, is a set of form $E = f^{-1}(\{ 1 \})$ for some $f \in A, \Vert f \Vert = 1$.  Equivalently, $E$ is a peak set iff there exists
 $g \in A$ with
$|g|_{|E^c} < \Vert g \Vert = 1 = g_{|E}$. A $p$-set is an
intersection of peak sets.   Hay defined a peak projection for a
unital subspace $A$ of a $C^*$-algebra $B$ to be a closed projection
in $B^{**}$, such that there exists an $a \in {\rm Ball}(A)$ with $a
q = q$ and
 satisfying any one of a long list
of equivalent conditions; for example $\Vert a r \Vert < 1$ for
every closed projection $r$ in $B^{**}$ with $r \leq q^\perp$. If
$A$ is a unital operator algebra, then peak projections are also the
complements of support projections of r-ideals in $A$ of the form
$\overline{(1-z)A}$ for $z \in {\rm Ball}(A)$ (see Proposition 6.7
in \cite{BHN}). By  \cite[Remark 6.10 (ii)]{BHN}, the latter support
projections  are the right support projections $r(1-z)$ for
contractions $z \in A$ (this also follows from results in Section 2
below).  A $p$-projection is defined to be the infimum of a family
of peak projections, or equivalently a weak* limit of a decreasing
net of peak projections.

If $A$ has a cai, then a {\em state} of $A$ is a functional $\varphi
\in {\rm Ball}(A^*)$ with $\varphi(e_t) \to 1$, for some (or every)
cai $(e_t)$ for $A$.  We write $S(A)$ for the space of states.  We
write $Q(A)$ for the quasistate space $\{ t \varphi : t \in [0,1],
\varphi \in S(A) \}$.  States extend uniquely to states on the unitization $A^1$
(see \cite[2.1.19]{BLM}).
We will sometimes use $C^*$-algebras
generated by an operator algebra $A$. If $C^*(A)$ is such a
$C^*$-algebra, then it is known that any bounded approximate identity
(bai) for $A$ is a bai for
$C^*(A)$, and hence states of $A$ are precisely the restrictions to
$A$ of states on $C^*(A)$ (see \cite[2.1.19]{BLM}). We will often use the
{\em numerical range} of an operator (see e.g.\ \cite{BD}), as
opposed to its spectrum.  This distinction is important: for
example, for an operator $T$,
 having spectrum $\{ 0 \}$ or contained in $[0,1]$ tells one very little,
whereas having numerical range in these sets gives $T = 0$ in the
first case, and $0 \leq T \leq I$ in the second. Of course the
(closed) numerical range of an operator contains its spectrum.

 For an operator algebra $A$, and $x \in A$, we define oa$(x)$ to be
the closed subalgebra of $A$ generated by $x$.
We define the {\em left} (resp.\ {\em  right})
{\em support projection} of $x \in A$ to be the smallest projection $p \in A^{**}$
such that $p x = x$ (resp.\ $x p = x$), if such a projection exists
(it always exists if $A$ has a cai).  If the left
and right support projections exist, and are equal,
 then we call it the {\em support projection}
written $s(x)$.

\medskip

We end this section with some quick consequences of Theorem \ref{read}:

Theorem  \ref{read} answers several questions posed in
\cite{Hay,BHN,ABS}. For example, it solves the biggest open problem
in Hay's thesis \cite{Hayth,Hay}.  This problem concerns
noncommutative peak sets, and the first part of the following result
may be viewed as the noncommutative version of a fundamental theorem
of Glicksberg on which the theory of peak sets rests (see
e.g.\ Theorem II.12.7 and II.12.5 in \cite{UA}).  Thus the
result puts the theory of noncommutative peak sets on a much firmer
foundation.

\begin{theorem} \label{corea1}  If $A$ is a unital operator algebra
and if $q$ is a closed projection in $A^{**}$, then
$q$ is a $p$-projection, and indeed is a strong limit of a
decreasing net of peak projections for $A$.

The r-ideals in a unital operator algebra $A$
are precisely the right ideals which are the
closure of the union of an increasing net of
right ideals of the form $\overline{(1-z)A}$ for $z \in {\rm Ball}(A)$.
\end{theorem}

\begin{proof}   The first statement was reduced, in the first two pages
of \cite[Section 6]{BHN},
to the existence in any operator algebra
with a cai, of a bai $(e_t)$ with $\Vert 1 - e_t \Vert \leq 1$
for all $t$.
The latter follows from Theorem \ref{read}.

For the second statement, we use the first statement, together
with the fact mentioned earlier that
peak projections are the complements of
support projections of r-ideals in $A$ of the
given form $\overline{(1-z)A}$ for $z \in {\rm Ball}(A)$.
It was shown in \cite[Proposition 6.8]{BHN}
that such  $\overline{(1-z)A}$ is an r-ideal (this also
follows from Lemma \ref{hasc} below).
 Moreover, the ordering of
open projections in $A^{**}$
corresponds to the inclusion of the matching r-ideals.
Hence, by the correspondence between r-ideals and open/closed
projections, closures of sums of r-ideals corresponds to
infs of closed projections (or sups of the complementary
open  projections).   More precisely, suppose that $(e_i)$ is a family of open projections
 corresponding to r-ideals $J_i$ in a
(possibly nonunital) operator algebra $A$.  Then $J$, the
closure of the span of the $J_i$, is known (and is easily seen) to be
an r-ideal, and its matching open projection $r$
equals $e = \vee_i \, e_i$.  Indeed $e \leq r$ clearly (since $J_i \subset J$).
Conversely, if $a \in J_i$ then $e_i a = a$, so that $e a = a$.
Hence $a = ea$ for any $a \in J$, so that
$r \leq e$, and $r = e$.

Putting this all together, any r-ideal is the closure of the
union of an increasing net of r-ideals of the
given form.
\end{proof}

{\bf Remark.}  In particular,
 every nonzero r-ideal in a unital operator algebra $A$ is what we called
1-{\em regular} in \cite{ABS}: that is it contains $(1-y)A$ for some $y \in {\rm Ball}(A)
 \setminus \{ 1 \}$.  This was stated as a question in that paper.

\begin{corollary} \label{rd0}  If $A$ is a nonunital operator algebra
with cai, and $x \in A^1 \setminus A$, then
there are always more than one  closest point in  $A$ to $x$.  That is,
$A$ is {\em never} a Chebychev subspace of $A^1$. \end{corollary}

\begin{proof}   The existence of nonzero $x \in A$
with $\Vert 1 - x \Vert = 1$ is saying that there are always more than
one  closest point in  $A$ to $1$.  If $a + \lambda 1 \in A^1$,
for $a \in A, \lambda \neq 0$, and  if $\Vert 1 - x \Vert = 1$ for $x \in A \setminus \{ 0 \}$, then
$\Vert a + \lambda 1 - (a + \lambda x) \Vert  = |\lambda| \Vert 1 - x \Vert = |\lambda|
= \Vert a + \lambda 1 - a  \Vert \leq ||a + \lambda 1 - b ||$
for all $b \in A$, using \cite[Lemma 2.1.12]{BLM}.  \end{proof}

See e.g.\ \cite{Ped2} for more information on Chebychev subspaces
 of operator algebras.

\begin{corollary} \label{rd1} Every r-ideal in an operator algebra, has a
left cai $(e_t)$ with
$\Vert 1 - 2 e_t \Vert \leq 1$ for all $t$, and $e_s e_t \to e_s$ with $t$ for any fixed $s$.
\end{corollary}

\begin{proof}  If $J$ is an r-ideal, and if $D$ is the matching HSA, then
by Theorem \ref{read}, $D$ has a  cai $(e_t)$ with
$\Vert 1 - 2 e_t \Vert \leq 1$.  Since $J = DA$, as explained in the
introduction, the result follows.
\end{proof}

\begin{corollary} \label{rd5}  If $J$ is a closed two-sided ideal in an operator algebra
$A$, and if $J$ has a cai, then $J$ has a cai $(e_t)$ with $\Vert 1 - 2 e_t \Vert \leq
1$ for all $t$, which is also quasicentral (that is,
$e_t a - a e_t \to 0$ for all $a \in A$).
\end{corollary}

\begin{proof}  Let  $(e_t)$ be a cai for $J$ with $\Vert 1 - 2 e_t \Vert \leq
1$ for all $t$ (see Theorem \ref{read}).  The weak* limit
$q$ of $(e_t)$ is well known to be a central projection in $A^{**}$, and
so $e_t a - a e_t \to 0$
weakly for all $a \in A$.  A routine argument using Mazur's theorem
shows that convex combinations of the $e_t$ comprise the desired cai,
and they will still have the property of being in the
convex set $\frac{1}{2} {\mathfrak F}_A$ defined earlier.  \end{proof}

\begin{corollary} \label{r6}  If $A$ is  an operator algebra with a  countable cai $(f_n)$, then $A$ has a countable cai in $\frac{1}{2} {\mathfrak F}_A$.
\end{corollary}

\begin{proof}
By Theorem \ref{read},
$A$ has a cai $(e_t)$ in $\frac{1}{2} {\mathfrak F}_A$.  Choosing $t_n$
with  $\Vert f_n e_{t_n} - f_n \Vert \vee \Vert e_{t_n} f_n  - f_n \Vert < 2^{-n}$,
it is easy to
see that $(e_{t_n})$ is a countable cai in $\frac{1}{2} {\mathfrak F}_A$.
 \end{proof}

\section{Consequences involving ${\mathfrak F}_A$}

\begin{lemma}   \label{hasc}
If $x \in {\mathfrak F}_A$, with $x \neq 0$,  then the operator
algebra {\rm oa}$(x)$ has a cai. Indeed, the operator algebra {\rm
oa}$(x)$ has a sequential cai belonging to $\frac{1}{2} {\mathfrak
F}_A$, consisting of elements $u_n=(\frac{x}{2})^{1/n}$, the $n$th
roots being suitably defined below.
 \end{lemma}

\begin{proof}  We will give two proofs  of the fact that  {\rm oa}$(x)$ has a cai,
 since both will be needed later.
The operator algebra oa$(x)$ is an ideal
in $C$, its unitization, which is the closed algebra generated by $1$
and $x$.    Indeed the closure of $x C$ is oa$(x)$.  Note
too that oa$(x)$ has a bai $(e_n)$
where $e_n = 1 - \frac{1}{n} \sum_{k=1}^n \,
(1-x)^k$, since $$\frac{x}{n} \sum_{k=1}^n \,
(1-x)^k = \frac{1}{n} (1-(1-x)) \sum_{k=1}^n \,
(1-x)^k = \frac{1}{n} (1 - (1-x)^{n+1}) \to 0 $$
with $n$.  Also, $\Vert \frac{1}{n} \sum_{k=1}^n \,
(1-x)^k \Vert \leq 1$.  By \cite[Theorem 6.1]{BHN}, oa$(x)$ has a
cai (the argument is that any weak* limit point $p$ of $(e_n)$ in
$A^{**}$ has to be the identity for oa$(x)^{**}$, hence is
idempotent.  Since $\Vert 1-p\Vert \leq 1$, we see that
$1-p$ and hence $p = 1_{{\rm oa}(x)^{**}}$ are projections.  So
{\rm oa}$(x)$ has a cai by a well known principle
stated in the introduction (see also \cite[Proposition 2.5.8]{BLM}).

The second proof will be presented after Proposition \ref{ob},
and it will include the extra information about the sequential cai $(u_n)$.
 \end{proof}

The following fact about the `disk algebra functional calculus' arising
 from von Neumann's inequality, is well known:

\begin{lemma}   \label{vnfc}
If $f, g \in A(\Ddb)$, with $\Vert g \Vert_{A(\Ddb)} \leq 1$,
and if $T \in B(H)$ is a contractive operator,
then $f(g(T)) = (f \circ g)(T)$.
 \end{lemma}

\begin{proposition}  \label{ob}   The sets  ${\mathfrak F}_A$ and $\frac{1}{2}
{\mathfrak F}_A$ are closed under taking roots.  That is, for $0 < r
\leq 1$ and $x\in   {\mathfrak F}_A$ (resp.\  $x\in \frac{1}{2}
{\mathfrak F}_A$),  a suitably defined $\rm{r}$th power $x^r$ is in
${\mathfrak F}_A$ (resp.\ $x^r \in \frac{1}{2} {\mathfrak F}_A$),
and $x^r\in {\rm oa}(x)$, and $(x^{r})^{\frac{1}{r}}=x$.
\end{proposition}

\begin{proof}   If $y\in A^1$ with $\Vert y \Vert \le 1$, then the disk algebra
functional calculus is a contractive algebra homomorphism $\theta$ from
$A(\Ddb)$ to $A^1$ with $\theta(1)=1$ and $\theta(z)=y$. If $r > 0$ then there is
a unique analytic branch of $f(z)=(1-z)^r$ defined on $\Ddb$
such that $f(0)=1$.  For  $x\in{\mathfrak F}_A$ set $y=1-x$.  Applying the functional
calculus for this value of $y$, $\theta(f)$ will be our suitable $r$th power of $x$.
The image $\theta(f)$ is a norm limit of polynomials
 $p_n(y)$, such that $p_n(z)$ converges
 uniformly to $(1-z)^r$ on the unit disk.
In particular the values at $z=1$ must tend to zero, and so
we may assume that $p_n(1) = 0$.  That is,
$x^{r}=\theta(f)$ is a norm limit of polynomials $q_n(x)$ with
$q_n(0)=0$, and these are in ${\rm oa}(x)$.   Hence   $x^{r}\in {\rm
oa}(x)$. Indeed for $0 < r \leq 1$, the binomial expansion of
$1-(1-z)^r$ is an absolutely convergent sum $\sum_{n=1}^\infty \,
a_n z^n$ with $a_n\ge 0$ and $\sum_{n=1}^\infty \, a_n=1$. Therefore
$1-x^r$ is in the closed convex hull of the powers $(1-x)^n$, so
$\Vert 1-x^r \Vert \le 1$ and
 $x^r\in{\mathfrak F}_A$.    A routine application of Lemma \ref{vnfc},
with $g(z) = 1-(1-z)^r$ and the $f$ there equal to $1 - (1-z)^{\frac{1}{r}}$,
 yields $(x^{r})^{\frac{1}{r}}=x$.

Suppose that $x \in \frac{1}{2} {\mathfrak F}_A$.  It is a pleasant exercise
in complex numbers that $\frac{1}{2} {\mathfrak F}_{\Cdb}$ is closed under
taking roots.  Equivalently, $|1 - 2(\frac{1-z}{2})^r| \leq 1$
for $0 < r \leq 1$ and $|z| \leq 1$.  Replacing $z$ by $1-2x$, that is
 by applying the functional calculus arising
 from von Neumann's inequality in a routine way, we have
 $\Vert 1 - 2 x^r \Vert \leq 1$.
\end{proof}

{\em Conclusion of proof of Lemma {\rm \ref{hasc}:}}  Suppose that
$x \in {\mathfrak F}_A$. It is not hard to see that $z^{\frac{1}{n}}
\; z \to z$  uniformly on the closed disk of radius $1$ center $1$.
Writing $y=1-x$, and
applying the functional calculus, we find that $\nm{x^{1/n}x-x}\to 0$.
The elements $u_n=x^{1/n}$ satisfy $u_nx=xu_n\to x$, and
so they are a bai for ${\rm oa}(x)$.
If $x\in \frac{1}{2} {\mathfrak F}_A$ then $u_n\in \frac{1}{2} {\mathfrak F}_A$,
and so $(u_n)$ is a cai.
 $\Box$

\begin{theorem} \label{newr}    For $0<\rho< \frac{\pi}{2}$
let $W_\rho$ be
the wedge-shaped region containing
the real interval $[0,1]$ consisting of numbers $re^{i \theta}$ with
argument $\theta$ such that $|\theta| < \rho$, which are also inside
the circle $|\frac{1}{2} - z | \leq \frac{1}{2}$.

An operator algebra $A$ with cai, has
a cai $(e_t)$ in $\frac{1}{2} {\mathfrak F}_A$,
 with the spectrum and numerical range of
$e_t$ contained in $W_\rho$.   In fact this can be
done with $\rho \to 0$ as $t$ runs over its directed set.
 \end{theorem}

\begin{proof}
If $x \in \frac{1}{2} {\mathfrak F}_A$
 then $x^{\frac{1}{k}}$ is in oa$(x)$ and in $\frac{1}{2} {\mathfrak F}_A$,
 by Proposition
 \ref{ob}, and it clearly has spectrum
contained inside a `wedge-shaped region' of the type described; and
the spectrum in $A$ is smaller.  The numerical range of $x$ is also
in this wedge, for example from a result of Macaev and Palant
\cite{MP} (see also e.g.\ \cite[Corollary 7.1.13]{sect}), stating that the
numerical range of a $k$th root of an
operator whose numerical range avoids the negative real axis, lies
in the appropriate `wedge' or sector centered on the positive real axis
of angle $\frac{\pi}{k}$.  It is also clearly inside the desired circle.

By Theorem \ref{read}, there is a cai $(u_t)$ in $\frac{1}{2} {\mathfrak F}_A$.
Let $v_{t,n} = u_t^{\frac{1}{n}}$ for $n \in \Ndb$.  If $b \in A$ then
using $(a_k)$ as we did in the proof of Proposition \ref{ob},
$$\Vert b - v_{t,n} \, b \Vert
= \Vert \sum_{k = 1}^\infty \, a_k (1-u_t)^k b\Vert \leq  (\sum_{k = 1}^\infty
 \, a_k)
\Vert (1-u_t) b \Vert = \Vert (1-u_t) b \Vert \to 0 ,$$
with $t$, for fixed $n$.  Similarly, $\Vert b - b v_{t,k}  \Vert \to 0$.
Thus $(v_{t,k})$ is also a cai in $A$.  By the last paragraph we can ensure
it has numerical range in the appropriate `wedges', and that
these wedges shrink to the interval $[0,1]$ with $(t,k)$.
\end{proof}

An operator with numerical range contained in $[0,1] \times
[-\epsilon , \epsilon]$, in fact is near to a positive operator.
Indeed Re$(x) = \frac{x + x^*}{2} \geq 0$ (since $\varphi(\frac{x + x^*}{2})
= {\rm Re}(\varphi(x)) \in [0,1]$ for states $\varphi$),
 and $\Vert x - {\rm Re}(x) \Vert
= \Vert {\rm Im}(x) \Vert \leq  \epsilon$ (since
${\rm Im}(x)$ is hermitian, so its norm is a supremum of quantities
$|\varphi(\frac{x - x^*}{2})|  =
|{\rm Im}(\varphi(x))| \leq \epsilon$).
It thus follows from Theorem \ref{newr} that {\rm any operator algebra with cai has a cai that gets
arbitrarily close to being positive}.  In fact this is not the deep
thing (the latter also follows by routine convexity methods of  \cite{AR,BAIC,Sm,MSIB}).
 What seems deep here is the position of the numerical range (being accretive
and sectorial, etc).

\begin{lemma}   \label{supp}   For any operator algebra $A$,
if $x \in {\mathfrak F}_A$, with $x \neq 0$,
then the left support projection of $x$
 equals the right support projection.
If $A \subset B(H)$ via a representation $\pi$, for a Hilbert space
$H$, such that the unique weak* continuous extension $\tilde{\pi} :
A^{**} \to B(H)$ is (completely) isometric, then $s(x)$ also may be
identified with the smallest projection $p$ on $H$ such that $p x =
x$ (and $x p = x$). That is, $s(x)H  = \overline{{\rm Ran}(x)} =
{\rm Ker}(x)^\perp$. Also,  $s(x)$ is an open projection in $A^{**}$
in the sense of \cite{BHN}. If $A$ is a subalgebra of a
$C^*$-algebra $B$ then $s(x)$ is open in $B^{**}$ in the sense of
Akemann {\rm \cite{Ake,Ake2}}.
 \end{lemma}

\begin{proof}  Viewing oa$(x) \subset A$, the identity of oa$(x)^{**}$
corresponds to a projection $e \in A^{**}$ with $e x = x e = x$.
If $A$ is represented on $H$ as described,
suppose that $p x = x$.  Then $p e_n = e_n$, where $(e_n)$ is the usual  bai
of oa$(x)$ from Lemma \ref{hasc}, so that in the weak* limit
we have $p e = e$ and $e \leq p$.  Similarly, $e \leq p$ if
$x p = x$.  So $e = s(x)$.  The equalities for $s(x)H$ are now routine.

This projection $e$, being the identity
of oa$(x)^{**}$, is open in the sense of \cite{BHN}.  The last statement
of the proof follows from e.g.\  \cite[Theorem 2.4]{BHN}.
 \end{proof}

\begin{corollary}  \label{supp2}  For any operator algebra $A$, if $x \in {\mathfrak F}_A$, with $x \neq 0$,  then
the closure of $xA$ is an r-ideal in $A$ and $s(x)$ is the support projection
of this r-ideal.  We have $\overline{xA} = s(x) A^{**} \cap A$.  Also,
$\overline{xAx}$ is the HSA matching $\overline{xA}$, and $x \in \overline{xAx}$.
\end{corollary}

\begin{proof}   The first assertion follows for example from
Lemma \ref{hasc}: any cai for oa$(x)$ serves as a left cai
for the closure of $xA$.  The second assertion follows  from
this, since the weak* limit of this left cai is $s(x)$.
Clearly $\overline{xA} \subset s(x) A^{**} \cap A$, and
since $(e_n)$ in the proofs above converges weak* to $s(x)$,
if $a \in s(x) A^{**} \cap A$ we have $e_n a \to a$ weakly.
By Mazur's theorem, a convex combination converges in norm, so
$a \in \overline{xA}$.

For the last assertion notice that by the argument in the first line of
this proof, $\overline{xAx}$ has a cai, and so it is a HSA.
It is the HSA matching $\overline{xA}$
by the correspondence described in the introduction,
since $\overline{xAxA} = \overline{xA}$.
The latter follows because $x \in \overline{xAx}$,
which in turn follows easily from Lemma \ref{hasc}.
 \end{proof}

\begin{corollary}  \label{perm}  If $A$ is a closed subalgebra of an
operator algebra $B$, and $x \in {\mathfrak F}_A$, then
the support projection of $x$ computed in $A^{**}$ is
the same, via the canonical embedding $A^{**} \cong A^{\perp \perp}
\subset B^{**}$, as the support projection of $x$ computed in $B^{**}$.
\end{corollary}

\begin{proof}   This is obvious given the formula  $s(x) = {\rm w^*}\lim_n \, e_n$ above.
 \end{proof}

\begin{corollary}  \label{perm2}  If $A$ is a closed subalgebra of a $C^*$-algebra
$B$, and $x \in {\mathfrak F}_A$, then $s(x)$ is the
support projection of $x^*x$ in $B^{**}$.   Indeed
$s(x) = s(x^* x) = s(x x^*) = s(x^*)$, where the latter three
support projections are with respect to $B$.
\end{corollary}

\begin{proof}
We have $x^*x s(x) = x^*x$, so $s(x) \geq s(x^*x)$.  Conversely,
if $p$ is a projection in $B$ with $x^*x p  = x^*x$, then
$(1-p) x^*x (1-p) = 0$, so that $x = x p$, and so $s(x) \leq p$
(using Corollary \ref{perm}).   Thus $s(x) \leq s(x^*x)$.  So $s(x) = s(x^* x)$
and the other equalities are similar, or now obvious.
 \end{proof}

\begin{lemma}   \label{su}  Let $A$ be an operator algebra with cai.
If $x \in {\mathfrak F}_A$,
then for any state $\varphi$ of $A$, $\varphi(x) = 0$ iff $\varphi(s(x)) = 0$.
 \end{lemma}

\begin{proof}  Let $B = C^*(A)$, then as we said in the
Introduction, states $\varphi$ on $A$ are precisely the restrictions
of states on $B$.  Continuing to write $\varphi$ for its
canonical extension to $A^{**}$, if $\varphi(s(x)) = 0$ then by Cauchy-Schwarz,
$$|\varphi(x)| = |\varphi(s(x) x)| \leq \varphi(s(x))^{\frac{1}{2}}
\varphi(x^* x)^{\frac{1}{2}} = 0 .$$
Conversely, if $\varphi(x) = 0$ then $\varphi(x^* x) \leq \varphi(x + x^*) = 0$,
since $x \in {\mathfrak F}_A$.
By Cauchy-Schwarz, $\varphi(ax) = 0$ for all $a \in A$.  Since any bai
for oa$(x)$ converges to $s(x)$ weak* we have $\varphi(s(x)) = 0$.
\end{proof}

\begin{lemma}   \label{supp3}  If $x, y \in {\mathfrak F}_A$, for
any operator algebra $A$, then $\overline{xA} \subset \overline{yA}$
iff $s(x) \leq s(y)$.  If $A$ has a cai and $x  \in {\mathfrak F}_A$,
then the following are equivalent:
  \begin{itemize}
\item [(i)] $\overline{xA} = A$.
\item [(ii)]  $\overline{xAx} = A$.
\item [(iii)]  $s(x) = 1_{A^{**}}$.
\item [(iv)]  $\varphi(x) \neq 0$ for every state $\varphi$ of $A$.
\item [(v)] ${\rm Re}(x)$ is strictly positive (that is,
$\varphi({\rm Re}(x)) > 0$ for every state $\varphi$ of $C^*(A)$).
 \end{itemize} \end{lemma}
\begin{proof}   Since $\overline{xA} = s(x) A^{**} \cap A$,
it is clear that
if $s(x) = 1$ then $\overline{xA} = A$; and also that
  $\overline{xA} \subset \overline{yA}$
if $s(x) \leq s(y)$.   Conversely, if $\overline{xA} \subset \overline{yA}$,
then $x \in \overline{yA} = s(y) A^{**} \cap A$.
We have $s(y) x = x$, so that $s(x) \leq s(y)$ by definition of $s(x)$.

(i) $\Leftrightarrow$ (iii) \ Corollary \ref{supp2}.

(i) $\Leftrightarrow$ (ii) \  Follows by the bijective correspondence between
r-ideals and HSA's, and Corollary  \ref{supp2}.

(iii) $\Rightarrow$ (iv) \ Obvious by the last Lemma.

 (iv)  $\Rightarrow$ (iii) \  If $s(x) \neq 1$ choose a state $\varphi$ on
$A$ (or equivalently on $C^*(A)$) with $\varphi(1-s(x)) = 1$.  Then
$\varphi(s(x)) = 0$, and so $\varphi(x) = 0$
by Lemma \ref{su}.

 (v)  $\Rightarrow$ (iv) \ Follows since Re$(\varphi(x)) = \varphi({\rm Re}(x))$.

(iv)  $\Rightarrow$ (v) \  If Re$(\varphi(x)) = \varphi({\rm Re}(x))
= 0$, then because $|1-\varphi(x)|\leq 1$, we must have $\varphi(x) = 0$.
\end{proof}

{\bf Remark.}  It is easy to see that in the last result we can replace $A$ by
any $C^*$-algebra $C^*(A)$ generated by $A$.  Thus for example
$\overline{x A} = A$ iff  $\overline{x C^*(A)} = C^*(A)$.

\bigskip

An element in ${\mathfrak F}_A$ with
${\rm Re}(x)$ strictly  positive, and hence
satisfying the equivalent conditions in
the last result, will be called {\em strictly real positive}.
Note that roots $x^{\frac{1}{k}}$ of a strictly real  positive
$x$ are strictly real positive,
and they become as close  as we like to a positive operator, as
$k \to \infty$.

\begin{proposition}  \label{srp}  If $x \in {\mathfrak F}_A$
is strictly real positive, then  $pxp$ is invertible in $pAp$ for every
projection  $p \in A$.
\end{proposition} \begin{proof}
The state space $S(pAp)$ is weak* compact since $pAp$ has an identity, and the map $\varphi \mapsto \varphi(pxp)$
on $S(pAp)$ is continuous.  It is also never zero, as can be seen
using Lemma \ref{supp3}, since
for any $\varphi \in S(pAp)$, $\varphi(p \cdot p)$
extends to a nonzero positive functional on $C^*(A)$, so is a
nonzero multiple of a state on $A$, hence $\varphi(pxp)
\neq 0$.   Thus $|\varphi(pxp)|$ is bounded
away from $0$, so $pxp$ has numerical range, hence spectrum with respect to
$pAp$, excluding $0$.  \end{proof}

We will need the following  `Fredholm alternative' type result,  a `sharp
 form' of the Neumann lemma.

\begin{theorem}  \label{abr}  Let $T$ be an operator in
 $B(H)$ with $\Vert I - T \Vert \leq 1$.  Then  $T$ is not invertible if
and only if
$\Vert I - T \Vert = \Vert I - \frac{1}{2} T \Vert = 1$.
Also, $T$ is invertible iff $T$ is invertible in the closed
algebra generated by $I$ and $T$, and iff  {\rm oa}$(T)$ contains $I$.
Here $I = I_H$
of course.
 \end{theorem}

\begin{proof}   Since $\Vert I - T \Vert \leq 1$ implies
$\Vert I - \frac{1}{2} T \Vert \leq 1$ by convexity, the $(\Rightarrow)$
direction of the first `iff' is clear by the Neumann lemma.
 Conversely, if $\Vert I - T \Vert = \Vert I - \frac{1}{2} T \Vert = 1$, then by the parallelogram law
$$\Vert \frac{1}{2} T  \zeta \Vert^2 + \Vert (I - \frac{1}{2} T) \zeta \Vert^2 = \Vert \frac{1}{2}  \zeta \Vert^2 + \Vert \frac{1}{2} (I-T) \zeta \Vert^2 \leq 1, \qquad \zeta \in {\rm Ball}(H) .$$  Hence $I - \frac{1}{2} T$
approximately achieves its norm at some norm one vector $\zeta$ with
$\Vert T \zeta \Vert$ as close as we wish to  $0$.  Hence $T$ is not invertible, or else
$\Vert T \zeta \Vert \geq \Vert T^{-1} \Vert^{-1}$.

If oa$(T)$ contains $I_H$ then $\Vert I -R T \Vert < 1$ for some
$R$ in oa$(T)$, which by commutativity of oa$(T)$
and the Neumann lemma implies that $T$ is invertible
in oa$(T)$, and hence in $A$.  Conversely, if
$T$ is invertible in $A$ then by the above,
$\Vert I - T \Vert < 1$ or $\Vert I - \frac{1}{2} T \Vert
< 1$.  In the first case, the bai
$(e_n)$ for oa$(T)$ in Lemma \ref{hasc},
 converges in norm to $I$, so $I \in {\rm oa}(T)$.  The second case
 follows from the first by replacing $T$ with $T/2$.
\end{proof}

\begin{lemma}  \label{genh}  If $(J_i)$ is a family of r-ideals in
an operator algebra $A$, with matching family of HSA's $(D_i)$,
and if $J = \overline{\sum_i \, J_i}$  then
the HSA matching $J$
is the HSA $D$ generated by the $(D_i)$ (that is, the smallest HSA in $A$ containing all the
$D_i$).  Here `matching' means with respect to the correspondence between
r-ideals and HSA's described in the introduction).
\end{lemma} \begin{proof}   Let $D'$ be the  HSA generated by the $(D_i)$.
Since $J_i \subset J$ we have $D_i \subset D$, and so $D' \subset D$.
Conversely, since $D_i \subset D'$ we have $J_i \subset D'A$,
so that $J \subset D'A$.  Hence $D \subset D'$.
\end{proof}

An r-ideal (resp.\
 HSA, $\ell$-ideal) of the form $\overline{x A}$ (resp.\ $\overline{xAx}$,
$\overline{Ax}$), for $x \in {\mathfrak F}_A$, will be called
{\em peak-principal}.  We note that the peak-principal r-ideals
in a uniform algebra $A$ are precisely the $J_E = \{ f \in A :
f_{|_E} = 0 \}$ for a peak set $E$.  Thus peak-principal r-ideals
(or peak-principal HSA's) may be thought of as a noncommutative
variant of peak sets (see also \cite[p.\ 354]{BHN}).

 If $J$ is a peak-principal r-ideal, for example, then
for any $\epsilon > 0$,
we may write $J = \overline{x A}$ for some
$x \in \frac{1}{2} {\mathfrak F}_A$, where  the numerical range of $x$
is contained in the thin wedge $W_\epsilon$ from
Theorem \ref{newr}. This
is because
  $\overline{x A} = \overline{x^{\frac{1}{n}} A}$ for all $n \in \Ndb$
  (since
$x^{\frac{1}{n}} \in {\rm oa}(x)$ by
Proposition \ref{ob}).

  As pointed out in \cite[Section 4]{BHN}, there is a bijective correspondence
between r-ideals, and certain weak* closed faces in the quasistate space
$Q(A)$, for an approximately unital operator algebra $A$.
 In fact, there are simple arguments for what we will need:
If $J$ is an r-ideal with support
projection $p$, let $F_p = \{ \varphi \in Q(A) : \varphi(p) = 0 \}
= Q(A) \cap J^\perp$.  (The one direction of the last equality
follows from Cauchy-Schwarz, the other from the fact that a left cai of $J$
converges to $p$.)    Note that if $F_{p_2} \subset F_{p_1}$
then by extending states to $C^*(A)$ we get a similar inclusion
with respect to $C^*(A)$.  Since $p_1, p_2$ are open with respect to
$C^*(A)$, we obtain $p_1 \leq p_2$ by the $C^*$-algebra theory,
and so $J_1 \subset J_2$.  Thus
 $J \mapsto F_p$ is a one-to-one order injection.   It also takes
`closures of sums of r-ideals' to `intersections', since
$Q(A) \cap (\sum_i \, J_i)^{\perp} = \cap_i \, (Q(A) \cap J_i^{\perp})$.

\begin{proposition}  \label{prim}  Let $A$ be any operator algebra
(not necessarily with an identity or approximate identity).
Suppose that $(x_k)$ is a sequence in ${\mathfrak F}_A$,
and that $\alpha_k \in (0,1]$ add to $1$.  Then the closure of the sum of the r-ideals
$\overline{x_k A}$, is the r-ideal $\overline{zA}$,
where $z = \sum_{k = 1}^\infty \, \alpha_k \, x_k \in {\mathfrak F}_A$.   Similarly,
the HSA generated by all the $\overline{x_k A x_k}$
equals $\overline{zAz}$.
\end{proposition}

\begin{proof}  Since $x \in {\rm oa}(x)$, it is easy to see
that $\overline{x A} = \overline{x A^1}$.  Thus we may assume that
$A$ is unital if we want.   The statement to be proved corresponds to
the identity $\vee_k \, s(x_k)  = s(z)$.
We have  $$\{ \varphi \in Q(A) : \varphi(\sum_k \, \alpha_k x_k) = 0 \} =
\cap_k \, \{ \varphi \in Q(A) : \varphi(x_k) = 0 \} $$
(because $\varphi(\sum_k \, \alpha_k x_k) = 0$ implies that
$\sum_k \, \alpha_k \, {\rm Re} \,  \varphi(x_k)  = 0$; and the latter
implies that ${\rm Re} \,  \varphi(x_k)  = 0$ since $x_k$ is accretive,
and so $\varphi(x_k)  = 0$ because of the shape of the numerical range
of elements in ${\mathfrak F}_A$).
Hence, by Lemma \ref{su}, $F_{s(z)} =
\cap_k \, F_{s(x_k)}$, which implies, in the light of the discussion
above the proposition, that the closure of the sum of the r-ideals
$\overline{x_k A}$, is the r-ideal $\overline{zA}$.    So
$\vee_k \, s(x_k)  = s(z)$.    The HSA assertion follows from the
r-ideal assertion, by Lemma  \ref{genh}, and the last assertions
of Corollary \ref{supp2}.
\end{proof}

\begin{theorem} \label{corea11}  Let $A$ be any operator algebra
(not necessarily with an identity or approximate identity).
The r-ideals (resp.\
 HSA's) in  $A$, are precisely the closures of unions of an increasing net of
of ideals (resp.\
HSA's) of the form $\overline{x A}$ (resp.\ $\overline{xAx}$), for
$x \in {\mathfrak F}_A$.
\end{theorem}

\begin{proof}
The r-ideal case is done in Theorem \ref{corea1} if $A$ is unital.
If $A$ is not unital, and if $J$ is an r-ideal in $A$, then it is also an r-ideal in $A^1$.
Theorem \ref{corea1} gives that
$J$ is the closure  of an increasing unions of ideals
of the form $(1-z)A^1$, for $z \in {\rm Ball}(A^1)$.
If $z = \lambda 1 - x$ for
$x \in A,  \lambda \in \Cdb$, then $\lambda
= 1$ (or else there is a nonzero scalar $t = 1- \lambda$ with
$t 1 + x = 1-z \in (1-z)A^1 \subset J \subset A$, which
forces $1 = \frac{1}{t}((t 1 + x) - x)
\in A$).  So $(1-z)A^1 = x A^1$.    Since oa$(x)$ has a cai
by Lemma   \ref{hasc}, $x \in \overline{x \, {\rm oa}(x)} \subset \overline{x A}$.
  It follows that
 $x A^1$ has the same closure as $x A$.  Thus the closure
 of $(1-z)A^1$ is $\overline{xA}$.  This completes the proof of the
 statements concerning r-ideals.

We saw in Corollary  \ref{supp2} that
$\overline{x A x}$ is a HSA if $x \in {\mathfrak F}_A$,
with matching r-ideal $\overline{x A}$.
  Also, the closure $D$ of an increasing union of HSA's
$D_i$ is an HSA, and therefore it is the HSA generated by the
$(D_i)$.  Indeed, it clearly satisfies
$DAD \subset D$, and to see that it has a cai it is well known
that it is enough to show that if $x_1, \cdots, x_n \in D$ and $\epsilon > 0$
are given, then there exists $d \in {\rm Ball}(D)$ with
$\Vert x_k d - x_k \Vert \vee \Vert d x_k - x_k \Vert < \epsilon$.
Picking $j$ and $y_1, \cdots, y_n \in D_j$ with
$\Vert x_k - y_k \Vert < \epsilon/3$ for all $k = 1, \cdots, n$,
we have for a cai $(e_t)$ for $D_j$ that $$\Vert x_k e_t  - x_k \Vert \leq
\Vert (x_k - y_k) e_t \Vert + \Vert y_k e_t- y_k \Vert +
\Vert x_k - y_k \Vert < \epsilon$$
for all $k$, and  for some choice of $t$.  Similarly,
$\Vert d x_k - x_k \Vert < \epsilon$, as desired.

  Finally, if $D$ is a HSA  in  $A$, with matching r-ideal $J$,
express $J$ as the closure of an increasing unions of r-ideals $J_i$
of the form $\overline{x A}$, by the first paragraph.  Then by Lemma
\ref{genh}, $D$ is the HSA generated by, and is also the closure of,
the increasing net of HSA's $(D_i)$ matching the $(J_i)$. By
Corollary  \ref{supp2}, $D_i$ is of the desired form $\overline{x A
x}$.
 \end{proof}

Peak sets were defined in the introduction.
It is well known that a countable intersection of peak sets is
a peak set, and that a p-set which is a G$_\delta$ set is a
peak set, so that for uniform algebras on metrizable spaces
the p-sets are exactly the peak sets (using the fact that
$C(K)$ is separable if $K$ is metrizable).  Analogous results hold
in $C^*$-algebras: e.g.\  separable closed right ideals are all
of the form $\overline{x A}$ for some $x \in A_+$.  Similar facts hold
in our context:

\begin{theorem} \label{sephs}  Let $A$ be any operator algebra
(not necessarily with an identity or approximate identity).
\begin{itemize}
\item [(1)]  Every separable r-ideal (resp.\
 HSA) in  $A$, is peak-principal (that is,
equal to  $\overline{x A}$ (resp.\ $\overline{xAx}$), for
some $x \in {\mathfrak F}_A$).
\item [(2)]  The closure of a countable sum
of peak-principal r-ideals (resp.\ The HSA generated by a countable
number of peak-principal HSA's)
is  peak-principal.
\end{itemize}
\end{theorem}

\begin{proof}
(2) \ By Proposition \ref{prim}, the closure of
$\sum_{n=1}^\infty \, \overline{x_n A}$ is $\overline{xA}$, where
$x = \sum_{n=1}^\infty \, \frac{x_n}{2^n}$.   The HSA assertion then
follows from this, as in the proofs of Proposition \ref{prim}
and Theorem \ref{corea11}.

(1) \ Let $D$ be a separable r-ideal (resp.\ HSA) in  $A$.
By Theorem \ref{corea11}, $D$ is the closure of a union of r-ideals (resp.\ HSA's)
 of the form $\overline{w A}$  (resp.\  $\overline{w A w}$), for $w  \in {\mathfrak F}_A$,
and we can clearly assume that there are a countable number of these.
 Now apply (2) (in the HSA case recall that the
HSA generated by an increasing net of HSA's is the closure of the union of these
HSA's).
 \end{proof}

{\bf Remark.}  The above considerations gives an `algorithm' for
building useful cai in r-ideals, $\ell$-ideals, or HSA's (and hence
in any approximately unital operator algebra). In the separable
case, we can just take $(y^{\frac{1}{k}})$ where $y = \frac{x}{2}$
for $x$ as in Theorem \ref{sephs}.
Indeed as we saw in the second
proof of Lemma \ref{hasc}, $y^{\frac{1}{k}} y \to y$.
Similarly in the
nonseparable case, any
r-ideal $J$, for example, in a unital operator algebra $A$ may
 be written as the closure of the union
of an increasing net of r-ideals $J_t = \overline{x_t A}$
for $x_t \in \frac{1}{2} {\mathfrak F}_A$, by Theorem
\ref{corea1}.   Then as before, $(x_t^{\frac{1}{k}})$
is a left cai for $J$.

\begin{corollary}    \label{ispk}  If $A$ is a separable
operator algebra, generating a $C^*$-algebra $B$, then
the open projections in $A^{\perp \perp}$  are precisely the
$s(x)$ for $x \in {\mathfrak F}_A$.
 \end{corollary}

\begin{proof}   This follows from Theorem  \ref{sephs} (1),
Lemma \ref{supp}, and Corollary \ref{supp2}.
\end{proof}

{\bf Remark.}  Of course if in the last result $A$ is also unital, then these
projections are precisely the
`perps' of peak projections for $A$, by a fact mentioned in the introduction.

\begin{corollary}    \label{csfa}  If $A$ is a separable
 operator algebra with cai, then there exists
an $x \in {\mathfrak F}_A$ with $A = \overline{xA} =
\overline{Ax} = \overline{x A x}$.

Any separable operator algebra with cai has a countable cai
consisting of mutually commuting elements, indeed
of form $(x^{\frac{1}{k}})$ for an $x \in \frac{1}{2} {\mathfrak F}_A$.
 \end{corollary}

\begin{proof}  The first part is immediate from
Theorem \ref{sephs}; if $A = \overline{x A x}$ then
this agrees with $\overline{x A}$ and $\overline{A x}$, since
for example  $\overline{x A x}
\subset \overline{x A} \subset A = \overline{x A x}$.   The second part
is clear
from  the fact that $x^{\frac{1}{k}} x \to x$ (see the second
  proof of Lemma \ref{hasc}).
 \end{proof}

\begin{theorem} \label{sap}  Let $A$ be any operator algebra
with cai.
The following are equivalent:
  \begin{itemize}
\item [(i)]  $A$ has a countable cai.
\item [(ii)]  $A$ has a strictly real positive element.
\item [(iii)]  There is an element $x$ in ${\mathfrak F}_A$ with $s(x) = 1_{A^{**}}$.
 \end{itemize}
\end{theorem}

\begin{proof}
 If $A$ has a strictly real positive element $x$ then $A =
\overline{xAx}$ by Lemma \ref{supp3}, and (a scaling of)
$(x^{\frac{1}{k}})$ is a  countable cai.

If $A$ has a countable cai $(f_n)$, then
$A$ has a countable cai $(e_n)$ in ${\mathfrak F}_A$ by
Corollary  \ref{r6}.
By Lemma \ref{supp3} and Theorem  \ref{sephs} (2),  $A = \overline{\cup_n \, e_{n} A} = \overline{z A}$
for a strictly real positive element  $z \in A$.

The equivalence of (ii) and (iii) comes from Lemma \ref{supp3}.
 \end{proof}

\begin{definition}  \label{320} If $A$ is an approximately
unital operator algebra, then we define a {\em peak projection} for $A$ to be the complement
of a support
projection $s(x)$, for an element $x \in {\mathfrak F}_A$.  A {\em $p$-projection} for $A$ is the infimum of a collection of peak projections for $A$.
\end{definition}

If $A$ is unital, this definition coincides with the ones discussed
in \cite{Hay, BHN}, as was pointed out in the introduction
(following from \cite[Remark 6.10]{BHN}).

\medskip

{\bf Remark.}  If $A$ is a nonunital $C^*$-algebra, then our
definition \ref{320} is connected, via Corollary  \ref{perm2}, to
the one in \cite{LNW}, but it is not the same. Indeed, the function
in $B = C_0((0,1))$ which is $0$ until $\frac{1}{2}$, and then makes
an inverted `vee' of height $1$, is in ${\mathfrak F}_B$, and the
corresponding peak projection is the characteristic function of
$[0,\frac{1}{2}]$.  However, the latter is not an essential support
projection in their sense.

\begin{corollary}    \label{supp5} For any approximately
unital operator algebra $A$, a projection $q \in A^{**}$
is the complement
of the support projection of an r-ideal iff $q$ is
the infimum of a collection of peak projections. These can be chosen to
be a decreasing net.  \end{corollary}   \begin{proof}
($\Rightarrow$) \ This follows  from Theorem \ref{corea11},
 and from the fact in the proof of Theorem \ref{corea1} that
the open projection corresponding to
the closure of a sum of r-ideals, is the supremum
of the open projections $p_t$ corresponding to each of these
r-ideals.  In our  case here
each $p_t = s(x_t)$ for some $x_t \in {\mathfrak F}_A$,
so that $p^\perp = \wedge_t \, s(x_t)^\perp$.

($\Leftarrow$)  \ This follows
from the fact that $s(x)^\perp$ is closed, and that the infimum
of closed projections in $A^{**}$ remains a
closed projection in $A^{**}$, as we said in the introduction.
   \end{proof}

Just as in the unital case \cite{Hay,BHN},
 one may write down several equivalent characterizations
of peak projections matching some of the characterizations in these papers.
We will not take the time to do this here since most of these become cumbersome
to state in the nonunital case.  We will mention a characterization
in terms of the tripotent $u(z) = {\rm w^*lim}_n \;  z (z^* z)^n$ considered by Edwards and Ruttiman \cite{ER4}:

\begin{proposition} \label{ruid}  If $A$ is any operator algebra
and if $x \in {\mathfrak F}_A$, set $z = 1-\frac{x}{2}$,
where $1$ is the identity of a $C^*$-algebra $B$ containing $A$.
Then $u(z)$ (computed
with respect to $B$) is a projection
and  $u(z) = s(x)^\perp = {\rm w^*lim}_n \; (z^* z)^n$.
\end{proposition}

\begin{proof}   We may assume that $A= B$,
and view $A \subset A^{**} \subset B(H)$.
Since $z = \frac{1 + (1-x)}{2}$, it is a contraction.
It is well known (and easy to see) that $(z^* z)^n \to P$
weak*, where $P$ is the projection onto Ker$(1-z^* z)$.  We claim that
Ker$(1-z^* z) = {\rm Ker}(x)$.  Indeed clearly ${\rm Ker}(x)
\subset {\rm Ker}(1-z^* z)$.  If $R$ is the contraction
$1-x$ and $\zeta \in H$ is a unit vector, then
$z^* z \zeta = \zeta$ implies that $(2 {\rm Re}(R) + R^* R) \zeta =
3 \zeta$, which implies that $\langle R \zeta , \zeta \rangle =
\langle R^* R  \zeta  , \zeta  \rangle = 1$.
Hence $\Vert R \zeta - \zeta \Vert^2 = 0$, so
$\zeta \in {\rm Ker}(x)$.    By
Lemma \ref{supp},  $s(x)^\perp$ is the projection onto
${\rm Ker}(x)$, namely $P$ above.  Of course $z P = (1 - x/2) P = P$.
Thus $(z^* z)^n \to s(x)^\perp$ weak*, and
 $z (z^* z)^n \to z P = P$ weak*.      \end{proof}

 \begin{corollary}    \label{ji} If $A$ is an approximately
unital operator algebra then every projection in $M(A)$ is
a strong limit of a decreasing net of peak projections, and is
also a strong limit of an increasing net of support projections
of elements of ${\mathfrak F}_A$.   If $A$ is separable then we do not need
to take limits here.
  \end{corollary} \begin{proof}
The first statements follow from  Corollary    \ref{supp5}, since every projection in $M(A)$ is
open and closed (see \cite[Proposition 5.1]{BHN}).
The last statement follows from this and Corollary \ref{ispk}.
 \end{proof}

During
the writing of the papers \cite{BHN,Hay}, we
had believed (on the basis of a proof that had a gap)
the following fact about compact projections $q$, hence
 about closed projections in the second dual
of a unital algebra: If $\{ u_i : i \in I \}$
is a collection of open projections whose supremum $\vee_{i \in I} \, u_i$
dominates
a compact projection
$q$, then $q \leq \vee_{i \in F} \, u_i$ for a finite set $F \subset I$.
This is true if $q = 1$, or under some strong commutativity hypotheses, but
is false in general (as may be seen by considering $A = \Kdb(\ell^2)$
(or its unitization if one prefers a unital algebra),
$q$  the projection onto $\Cdb e_1$, and $u_k$  the projection onto
Span$(\{ e_1 + e_2, e_2 + e_3 , \cdots , e_k + e_{k+1} \})$.
Then $\vee_{k} \, u_k = I$, but we do not have
$q \leq \vee_{k = 1}^n \, u_k$ for any finite $n$).

This incorrect statement was used only twice in those
papers, namely in \cite[Proposition 5.6]{Hay} and   \cite[Theorem 6.4]{BHN}.
Fortunately both of these proofs can be fixed.
There is a very short direct proof for \cite[Proposition 5.6]{Hay}:
 note that the result is true
for peak projections since these are weak* limits of terms in $A$.
Every p-projection $q$ is a limit of a decreasing net
of peak projections $q_i$,
so $\varphi(q) = \lim_i \, \varphi(q_i) = 0$.
We can fix the gap in the first two lines
of one direction of the proof of \cite[Theorem 6.4]{BHN},
and at the same time  improve the result as follows:

\begin{theorem} \label{peakch} {\rm (Noncommutative Urysohn lemma for
 nonselfadjoint operator algebras) } \ Let $A$ be a
subalgebra of a unital $C^*$-algebra $B$, with $1_B \in A$, and let
$q \in B^{**}$ be a closed projection.  Then
$q \in A^{\perp \perp}$ iff for any open
projection $u \geq q$, and any $\epsilon > 0$,
there exists an $a \in {\rm Ball}(A)$ with
$a q = q$ and $\Vert a (1-u) \Vert < \epsilon$
and $\Vert (1-u)  a \Vert < \epsilon$.
Indeed this can be done with, in addition, $\Vert 1 - 2 a \Vert \leq 1$.
\end{theorem}

 \begin{proof}
($\Leftarrow$) \   As in \cite[Theorem 6.4]{BHN}.

($\Rightarrow$) \  Let $q  \in A^{\perp \perp}$, let $u$ be an open
projection with $u \geq q$, and let $\epsilon > 0$ be given.   Using
Theorem \ref{read}, let $(e_t)$ be a cai
with $\Vert 1 - 2 e_t \Vert \leq 1$,
for the HSA $q^\perp A^{**} q^\perp \cap A$
associated with $q$ as in the introduction.
Then $1 - e_t \to q$ weak*, and $(1 - e_t) q = q$.  We follow
the idea in the last seven lines
of the proof of \cite[Theorem 6.4]{BHN}:
By the noncommutative Urysohn lemma
\cite{Ake}, there is an $x \in B$ with
$q \leq x \leq u$.
Then  $(1 - e_t) (1-x) \to q (1-x) = 0$ weak*, and hence weakly
in $B$.  Similarly, $(1-x) (1 - e_t) \to  0$ weakly.
By a routine convexity argument in $B \oplus B$, given $\epsilon > 0$
there is a convex combination $a$ of the $1 - e_t$
such that $\Vert a (1-x) \Vert < \epsilon$ and
$\Vert  (1-x) a \Vert < \epsilon$.  Clearly
$\Vert 1 - a \Vert \leq 1$ and $a q = q$.  Therefore $\Vert a (1-u) \Vert =
\Vert a (1-x) (1-u)  \Vert < \epsilon$.
Similarly for $\Vert   (1-u) a \Vert  < \epsilon$.
The estimate $\Vert 1 - 2 a \Vert \leq 1$ follows
since $\Vert 1 - 2 (1 - e_t) \Vert  = \Vert 1 - 2 e_t \Vert  \leq 1$,
and this formula persists with convex combinations.
 \end{proof}

We now show that the elements $a$ in Theorem \ref{peakch}
constitute a left cai for the r-ideal associated with $q$,
with the net constituting the cai
indexed by the directed set of open projections $u \geq q$.

\begin{corollary} \label{peakchop}   Let $A$ be a
unital-subalgebra of $C^*$-algebra $B$ and let
$q \in A^{\perp \perp}$ be a closed projection associated
with an r-ideal $J$ in $A$.
Then an explicit left cai for $J$
is given by   $x_{(u,\epsilon)} = 1-a$, where $a$ is
an element which satisfies the
conclusions of the last theorem, and is associated
to an open projection $u \geq q$,
and a scalar $\epsilon > 0$.
This left cai is indexed by such pairs $(u,\epsilon)$,
that is, by the product of the directed set of open projections $u \geq q$,
and the set of $\epsilon > 0$.    This right cai is also in
$\frac{1}{2} {\mathfrak F}_A$; that is,
 $\Vert  1 - 2 x_{(u,\epsilon)} \Vert \leq 1$.
\end{corollary}

 \begin{proof}
By \cite[Theorem 2.9]{AP}, if $p_1, p_2$ are two
open projections dominating $q$, then there exists a third
open projection $p \geq q$ with
$\{ x \in B_+ : q \leq x \leq p \} \subset \{ x \in B_+ : q \leq x \leq p_k \}$.
By Lemma 2.7 in \cite{AP}, there is an increasing net in the first of these
sets with strong limit $q$, hence $q \leq p_k$ for $k = 1, 2$.  Thus
the set of open projections dominating $q$ is a directed set.

By the last theorem, for each open projection $u \geq q$,
and any $\epsilon > 0$,
there exists an $a \in {\rm Ball}(A)$ with
$a q = q$ and $\Vert a (1-u) \Vert < \epsilon$, and
 $\Vert 1 - 2 a \Vert \leq 1$.  Let $x_{(u,\epsilon)} = 1-a$.  We claim that $(x_{(u,\epsilon)})$
is a cai for $J$.  Certainly $x_{(u,\epsilon)} q = (1-a) q = q - q = 0$, so that
$x_{(u,\epsilon)} \in J$.  Also, $\Vert x_{(u,\epsilon)} \Vert \leq 1$, indeed
$\Vert 1 - 2 x_{(u,\epsilon)} \Vert \leq 1$.  If $b \in {\rm Ball}(J)$, we have
 $$\Vert x_{(u,\epsilon)} \, b - b \Vert = \Vert a b \Vert \leq
\Vert (a - au) b \Vert + \Vert au b \Vert \leq \Vert a(1 - u) \Vert +
 \Vert u b \Vert < \epsilon + \Vert u b \Vert .$$
Thus we need to show that $u b \to 0$ in norm with $u$, over the directed set of
open  projections $u \geq q$.  This is easy, for example in
Akemann's proof in \cite[Proposition 2.3]{Hay} one associates to an
increasing left cai $(a_t)$ for  $q^\perp B^{**} \cap B$, an open
projection $r_t$ with $q  \leq r_t \leq 2(1-a_t)$.
It follows that $\Vert r_t b \Vert \leq 2 \Vert (1-a_t) b \Vert \to 0$, since
$b \in J \subset q^\perp B^{**} \cap B$.  Thus if $\Vert r_t b \Vert < \epsilon$,
 then $$\Vert x_{(u,\epsilon)} b - b \Vert < \epsilon + \Vert u b \Vert
< \epsilon + \Vert r_t  b \Vert < 2 \epsilon ,$$
if $q \leq u \leq r_t$.
\end{proof}

\section{When $xA$ and $Ax$ are closed}

\begin{proposition}  \label{obc}  If $A$ has a cai but no identity,
and $x  \in {\mathfrak F}_A$ with $\overline{A x} = A$, then $xA \neq A$.
Hence for no strictly real positive $x \in A$ is $xA$ closed.
\end{proposition}

\begin{proof}  If $x = xy$ for some $y \in A$,  then
$e_t y = e_t \to y$ for the cai $(e_t)$,  so that $A$ has identity.
For the last statement use Lemma \ref{supp3}.  \end{proof}

We recall that
`well supported' operators are those operators $x$
 that have a `spectral gap'  (for $|x|$)
at $0$, that is $0$ is absent from, or is isolated in, the
spectrum (of $|x|$).   It is a well known principle in operator theory and
$C^*$-algebras, that
$x$ is a well supported operator (resp.\ element
if a $C^*$-algebra $A$) iff $x$ has closed range (resp.\
 $xA$ is closed), and this is equivalent to the existence of
an operator $y$ (resp.\ element $y \in A$) with $x y x = x$.  Such a $y$ is
called a {\em generalized inverse} or {\em pseudoinverse}.
See e.g.\ \cite[II.3.2.11]{Bla}, \cite{HM,HM2}.    With this in mind,
it is tempting to conjecture that for an operator algebra $A$,
a noninvertible element $x \in {\mathfrak F}_A$ has $0$ isolated in (or absent from)
its spectrum, iff $xA$ is closed, and iff there exists
$y \in A$ with $x y x = x$.
  However there are two issues that we have to deal with, for
$x \in {\mathfrak F}_A$.
First, we do not know if it is
 true that $xA$ is closed iff $Ax$ is closed.
What is true is
that $xA$ and $Ax$ is closed iff $xAx$ is closed.
Second, in a nonsemisimple setting,
 $0$ being an isolated point
in Sp$(x)$ need not imply that $xA$ is closed. Indeed there can
exist quasinilpotent operators without closed range.  Thus for
example suppose that $A$ is the radical operator algebra in Example
\ref{Vol}, with cai $(e_t) \subset {\mathfrak F}_A$. Then $e_t$ is
quasinilpotent (any character of $A$ certainly annihilates $e_t$, so
Sp$(e_t) = (0)$).  Hence  $0$ is an isolated point of its spectrum,
but $\overline{e_t A} = A$ since this algebra has no proper
r-ideals, and this differs from $e_t A$ by Proposition \ref{obc}. So
$e_t A$ is not closed.

With the above in mind, the following result may be the best possible:

\begin{theorem}   \label{ws}  For
an operator algebra $A$, if $x \in {\mathfrak F}_A$, the following
are equivalent:
  \begin{itemize}
\item [(i)] $x \, {\rm oa}(x)$ is closed.
\item [(ii)]  ${\rm oa}(x)$ is unital
(which implies that $x$ is invertible in ${\rm oa}(x)$).
\item [(iii)]  There exists $y \in {\rm oa}(x)$ with $xy x = x$.
\item [(iv)] $xAx$ is closed.
\item [(v)]  $xA$ and $Ax$ are closed.
\item [(vi)]  There exists $y \in A$  with $xy x = x$.
\end{itemize}
Also, the latter conditions imply
 \begin{itemize}
\item [(vii)]  $0$ is isolated in, or absent from, Sp$_A(x)$.
\end{itemize}
Finally, if further ${\rm oa}(x)$ is semisimple, then
conditions {\rm (i)--(vii)} are all equivalent.
\end{theorem}

 \begin{proof}  If $A$ is unital, $x \in {\mathfrak F}_A$, and $x$ is
invertible in $A$, then by Theorem \ref{abr} we have (ii), and
indeed in this case (i)--(vii) are all obvious.  So we can assume
that $x$ is not invertible in $A$.

That (i) implies  (ii) follows since in this case $x \in \overline{x \, {\rm oa}(x)} =
x \, {\rm oa}(x)$, and if $x = xy$ for $y \in {\rm oa}(x)$ then $y = 1_{{\rm oa}(x)}$.
That the first condition in (ii) implies the second is obvious (as in
Theorem \ref{abr}).
Now the equivalences (i)--(iii)
are obvious (some also follow from (iv)--(vi)), as is the fact that these
imply some of (iv)--(vi).

(iv)  $\Rightarrow$ (vi) \ Suppose that $xAx$ is closed.
 Now  $x = (x^{\frac{1}{3}})^3$, and
$$x^{\frac{1}{3}} \in {\rm oa}(x) = \overline{x \, {\rm oa}(x) \, x}
\subset \overline{x \, A \, x} = xAx ,$$
and so $x = x y x$ for some $y \in A$.

(vi)  $\Rightarrow$ (v) \ (vi) implies
that $xA = xy A$ is closed since $xy$ is idempotent.  Similarly
$Ax$ is closed.

(vi) $\Rightarrow$ (iv) \   $xAx = xy A yx$ is closed as in the last line.

(v)  $\Rightarrow$ (vi) \ If $xA$ and $Ax$ are closed,
then by a similar argument, $x^{\frac{1}{3}} \in xA$, and similarly
$x^{\frac{1}{3}} \in Ax$.  Hence
$x = x^{\frac{1}{3}} x^{\frac{1}{3}} x^{\frac{1}{3}}
\in xAx$.  Thus $x = x y x$ for some $y \in A$.

(vi) $\Rightarrow$ (vii) \ We may assume that $A \subset B(H)$, and
that $x$ is pseudo-invertible as an operator on $H$. Then $x(H) \subset xy(H) \subset
x(H)$, so that $x(H)$ is closed. Let $P$ be the projection onto $K =
{\rm Ker}(x)^\perp = x(H)$ (see Lemma \ref{supp}), an invariant subspace for
$x$.  Let $S$ be the restriction of $x$ to  $K$, then  $S$
is bicontinuous and  invertible.
If $K = H$ then $x$ is invertible,
and we discussed this case at the start of the proof.
If $K \neq H$ then since $x = P S P$, it follows that Sp$_{B(H)}(x)
= \{ 0 \} \cup {\rm Sp}_{B(K)}(S),$ which has $0$ as an isolated
point since $S$ is invertible.  If $0$ was not isolated in ${\rm
Sp}_{A}(x)$ then, by the topology of compact sets in the plane,
there is a sequence of boundary points in ${\rm Sp}_{A}(x)$
converging to $0$. Since $\partial {\rm Sp}_{A}(x) \subset {\rm
Sp}_{B(H)}(x)$,
 this is a contradiction.

(vi) $\Rightarrow$  (ii) \  We saw in the last paragraph or two that
(vi) implies that $x = 0 \oplus S$ for an invertible $S$. Note that
oa$(x) \cong {\rm oa}(S)$.  By Theorem \ref{abr}, ${\rm oa}(S)$ is
unital, and $S$ is invertible there. Thus the same is true for $x$.

(vii) $\Rightarrow$  (ii) \ Suppose that ${\rm oa}(x)$ is semisimple but nonunital,
and that $0$ is an isolated point in Sp$_A(x)$.
The latter is equivalent
(by the basic spectral result for
singly generated subalgebras, and
because the spectrum is contained in the disk $B(1,1)$),
to $0$ being isolated  in
 $K = {\rm Sp}_{{\rm oa}(x)}(x)$.
Consider $E$, the spectral projection of $x$ corresponding
to $\{ 0 \}$.  Namely, $E = f(x)$ where
$f$ is $1$ on a neighborhood of $0$, and is $0$
 on a neighborhood of $K \setminus \{0 \}$.
We have ${\rm Sp}(Ex) = {\rm Sp}((f z)(x)) =
(f z)({\rm Sp}(x)) = \{0 \}$.    By semisimplicity,
$Ex = 0$, and $(1-E) x = x$.  Since $(1-f)(0) = 0$ we
have $1 - E \in {\rm oa}(x)$.   Note that
if $g$ is $0$ on a neighborhood of $0$, and is $1/z$ on
a neighborhood of $K \setminus \{0 \}$, then $g z^2 - z$ is
zero on ${\rm Sp}(x)$, and so by semisimplicity we have
$x g(x) x = x$.   Since $g(0) = 0$
we have $g(x) \in {\rm oa}(x)$, and we deduce
from the above that ${\rm oa}(x)$ is unital, and (ii) holds.
 \end{proof}

{\bf Remark.}  The conditions in the theorem are not necessarily equivalent
under the assumption that $A$ is semisimple.  For example, if $A = B(L^2([0,1]))$,
suppose that $T \in A$ is any quasinilpotent operator with $T - I$ contractive.
Then  $T$ does not have closed range, for if it did
have closed range then as in the proof that (vi) $\Rightarrow$ (vii)
above, $T$ is of the form $0 \oplus S$ for an invertible $S$.  This
is impossible for a quasinilpotent (since if $0 \neq t \in {\rm Sp}(S)$ then
$t  \in {\rm Sp}(T) = (0)$, a contradiction).
It is easy to see that then (vi) fails, and so (i)--(v) fail too.
However $0$ is isolated in the spectrum of this quasinilpotent operator.

Also in this connection we remark that all $C^*$-algebras are semisimple,
yet $C^*$-algebras may have dense subalgebras consisting entirely of
nilpotents (hence quasinilpotents) \cite{Read0}.

\section{Operator algebras without HSA's}

In this section we study operator algebras $A$ without nontrivial HSA's or r-ideals.
By a trivial HSA (or r-ideal) of $A$ we mean of course $(0)$ or $A$.

\begin{theorem} \label{nor}
For a unital operator algebra $A$, the
  following are equivalent:
\begin{enumerate}
\item [{\rm (i)}]  $A$ has no nontrivial r-ideals (or equivalently, HSA's).
\item [{\rm (ii)}]  $a^n \to 0$ for all $a \in {\rm Ball}(A) \setminus \Cdb 1$.
\item [{\rm (iii)}] The spectral radius  $r(a) < \Vert a \Vert $  for all $a \in {\rm Ball}(A) \setminus \Cdb 1$.
\item [{\rm (iv)}] The numerical radius
 $\nu(a) < \Vert a \Vert $  for all $a \in {\rm Ball}(A) \setminus \Cdb 1$.
\item [{\rm (v)}]  $\Vert 1  + a \Vert < 2$ for all $a \in {\rm Ball}(A)\setminus \Cdb 1$.
\item [{\rm (vi)}]    ${\rm Ball}(A)\setminus \Cdb 1$ consists entirely of quasi-invertibles.
 \end{enumerate}
 If $A$ has a cai but no identity, then the
  following are equivalent:
 \begin{enumerate}
\item [{\rm (a)}]  $A$ has no nontrivial r-ideals (or equivalently, HSA's).
    \item [{\rm (b)}]
$A^1$ has one nontrivial r-ideal.
\item [{\rm (c)}]  ${\rm Re}(x)$ is strictly positive
for every $x \in {\mathfrak F}_A \setminus \{ 0 \}$.
\end{enumerate}
\end{theorem}

\begin{proof}  The equivalences (i)--(vi) are in \cite[Section 3]{ABS}, if
one uses the fact that every r-ideal contains what we called
a 1-regular ideal (defined before
Corollary \ref{rd0}), which is a consequence of Theorem \ref{read}.
The one direction of the equivalence of (a) and (b) is obvious.   For the other,
suppose that $x \in {\rm Ball}(A^1)$.  Then $\overline{(1-x) A}$
is an  r-ideal in $A$ by \cite[Proposition 3.1]{BHN} (or Theorem \ref{best}
below).
  So if $A$ has no nontrivial r-ideals
then either $(1-x)A = (0)$ or $\overline{(1-x) A} = A$.
In the first case, $x a = a$ for all $a \in A$,
which forces (if $A$ has a cai) $x = 1$ and
 $\overline{(1-x) A^1} = (0)$.  In the
second case: in the notation
of the proof of Theorem \ref{best} and the remark after
it, $\overline{(1-x) A}$ has bai $(e_n f_t)$,
which has weak* limit $p f$
if $f$ is the left identity of $A^{\perp \perp}$ and
$p$ is the weak* limit of $e_n$.  So $pf = f$,
which forces $p = f$ or $p = 1$.  If $p = f$
then $\overline{(1-x) A^1} = A$.  If $p = 1$ then $\overline{(1-x) A^1} =
A^1$.    By the remark after
Theorem \ref{corea1}, for example, $A^1$ has $A$ as its only  nontrivial r-ideal.

That (a) is equivalent to (c) follows easily from Lemma \ref{supp3},
 and
the fact that r-ideals are `sups' of peak-principal ones (see
Theorem \ref{corea11}).
\end{proof}

{\bf Remarks.}   1) \ Simple examples (the $2\times 2$ matrices
supported on the first row) shows that the equivalence with (b) in
the last result is not true without a cai.

\smallskip

2) \ If there exist nontrivial r-ideals in a unital
operator algebra $A$,  then there exist proper
maximal r-ideals in $A$.  This follows from \cite[Proposition 3.6]{ABS}
and the Remark after Theorem \ref{corea1}
above.

\begin{theorem}    \label{cc}  An approximately unital operator algebra
with no countable cai, has nontrivial r-ideals.
\end{theorem}

\begin{proof}  If $A$ has no countable cai then by Theorem \ref{sap}
there is no element in ${\mathfrak F}_A$ with $s(x) = 1$.
Thus for any nonzero $x \in {\mathfrak F}_A$, we have
$\overline{xA} \neq A$ by Lemma \ref{supp3}, and this is a nontrivial r-ideal.
 \end{proof}

\begin{proposition}    \label{cc2}  If a nonunital operator algebra
$A$ contains a nonzero $x \in {\mathfrak F}_A$ with $xAx$ closed,
or with $0$ isolated in ${\rm Sp}_A(x)$ and ${\rm oa}(x)$
semisimple, then $A$ has a nontrivial r-ideal.
\end{proposition}

\begin{proof}  By Theorem   \ref{ws}, under these conditions
 $x A$ is closed, and so $\overline{xA}
\neq A$ by Proposition \ref{obc}.   This is a nontrivial r-ideal.
\end{proof}

\begin{proposition}    \label{ccra}   A  nontrivial r-ideal
in the unitization of an approximately unital radical
operator algebra $A$,  is an  r-ideal in $A$.
\end{proposition}

\begin{proof} If $\lambda \neq 1$ then we claim that $z = \lambda 1 + a$
is quasi-invertible in $A^1$, so that
$(1-z) A^1 = A^1$,  for all $a \in A$.  Indeed we know that
$\frac{a}{1-\lambda}$ is quasi-invertible in $A$, and easy algebra
shows that its quasi-inverse gives rise to a quasi-inverse of
$\lambda 1 + a$.   \end{proof}

We now give several examples of
operator algebras with  cai, with only trivial r-ideals.

\subsection{Example. A unital two
dimensional example}  Consider the upper triangular $2 \times 2$ matrices
whose $1$-$2$ entry is the difference of the diagonal entries.

\subsection{Example. A  nonunital commutative semisimple example}
\label{uss}  Set   $A = R D R^{-1
}$, where
$D$ is
the diagonal copy of $c_0$   in $B(\ell^2)$, and $R$ is an invertible operator in
$B(\ell^2)$, such that the
commutant of $R^* R$ contains no nontrivial projections in
the diagonal copy of $\ell^\infty$ in $B(\ell^2)$. For example,
$R = I + S/2$ where $S$ is the backwards
shift.  Since $A$ is a
subalgebra of the compact operators, its second dual may be identified
with its $\sigma$-weak closure in $B(\ell^2)$.  Thus $A^{**}$
is unital, so that $A$ has cai.     In this case, there are no nontrivial projections
in $A^{**} = A^{\perp \perp} = R \overline{D}^{w*} R^{-1}$.  Indeed
any  projection $q$ in $R \overline{D}^{w*} R^{-1}$ corresponds to an idempotent,
hence projection, $p$ in $\overline{D}^{w*}$.  Any projection in $\overline{D}^{w*}$
 is a sequence of $0$'s and $1$'s.   That $q = q^*$ forces
$p$ to commute with $R^* R$, so that $p = 0$ or $p=1$. Thus $A$ has no
nontrivial r-ideals.

\subsection{Example.  A commutative radical example}
 \label{Vol}  Let $A$ be the norm
closed algebra generated by
the
Volterra operator $V$.  This is commutative, and so since $V$
is quasinilpotent we have that $A$ is radical.
By \cite[Corollary 5.11]{Dav} we have $\bar{A}^{w*} = V'$. As in the last example
this coincides with $A^{**}$, and since this is unital we see that $A$ has
a cai.
 By  \cite[Lemma 5.1]{Dav}, $V'$
contains no nontrivial projections, hence the same is true for $A^{**}$.
Thus $A$ has no HSA's or r-ideals.

\bigskip

In the next section we will continue looking at examples.

\section{An approximately unital radical operator algebra which is
an integral domain}  \label{cdv}

In this section we present an interesting commutative
approximately unital operator algebra, which happens to be radical
 and semiprime,
and in fact is an integral domain (so $ab = 0$ exactly when $a = 0$ or
 $b = 0$).
It is also an operator algebra whose ideal structure we can completely describe.
This is achieved by adapting the work of Domar \cite{domar}
on convolution algebras $L^1(\bR^+,\omega)$; he showed that certain conditions on
$\omega$ imply that $L^1(\bR^+,\omega)$ is a radical Banach algebra and that
all of its closed ideals are of a certain `standard' type
which we will describe below. A simplified exposition of Domar's result
can be found on p.\ 554 of Dales \cite{Dal}.
The operator algebras we produce will clearly have no nontrivial ideals
having approximate identities.  They will have some
features in common with Example \ref{Vol}, but in other ways they are very different.

For our purposes, a {\it radical weight} $\omega:[0,\infty)\to(0,
\infty)$ is a continuous function such that $\omega(0)=1$, $\omega(
s+t)\le\omega(s) \; \omega(t)$ for all $s,t\ge 0$,
 and $\omega(t)^{1/t}\to 0$ as $t\to\infty$.
The Banach spaces $L^p(\bR^+,\omega)$  ($1\le p<\infty$) consist of
 equivalence classes of measurable functions $f:\bR^+\to\bC$  such that
 ${\nm f}_p=(\int_0^\infty \,
|f(t)|^p \, \omega(t)^p \, dt)^{1/p}<\infty$. The space $L^1(\bR^+
,\omega)$ is a  Banach algebra when given the convolution multiplication (see \cite[Section 4.7]{Dal}).
  For each $\alp\ge 0$, there is the ``standard ideal'' $J_\alp\subset
L^1(\bR^+,\omega)$ consisting of functions supported on $[\alp,\infty)$, and this ideal is always norm closed.

We say that the radical weight $\omega$ satisfies Domar's criterion
if the function $\eta(t)=-\log\omega(t)$ is a convex function
on $(0,\infty)$,
and for some $\veps>0$ we have $\eta(t)/t^{1+\veps}\to\infty$ as $t\to\infty$. An obvious example of such a weight is $\omega(t)=e^{-t^2}$. Domar's theorem asserts that if the radical weight $\omega$ satisfies Domar's criterion, then the standard ideals are the only nonzero closed ideals in $L^1(\bR^+,\omega)$.
We will use this result to obtain radical operator algebras with interesting properties.
Let $\omega$ denote any radical weight.  The algebra $L^1(\bR^+,\omega)$ acts on
$L^2(\bR^+,\omega)$ by convolution, for one may readily check that the familiar inequality $\nm{f*g}_2\le {\nm f}_1\cdot {\nm g}_2$ still holds when the $L^p$ spaces are given their radical weighting according to $\omega$.
If we write $H$ for the Hilbert space $L^2(\bR^+,\omega)$ and   $M_f$ for the operator on $H$ with $M_f(g)=f*g$, then the norm closure of the operators $M_f\ (f\in  L^1(\bR^+,\omega))$, is an operator algebra $\cA=\cA(\omega)$.  Indeed,
 the set of operators $M_f$ is already a subalgebra of $B(H)$ since $M_f\cdot M_g=M_{f*g}$.

Clearly $\cA$ is a commutative operator algebra, and we claim that it has a dense set of
quasinilpotent elements consisting of operators $M_f$, so that $\cA$ is also radical.
Let $|\cdot|_1$ denote the {\it unweighted} $L^1$ norm, that is,
the usual norm on $L^1(\bR^+)$, and let $\Nm_1$ be
 the weighted norm on $L^1(\bR^+,\omega)$.
If $f\in L^1(\bR^+,\omega)$
 is supported on $[a,b]$, where $0<a<b$, then ${\nm f}_1\ge |f|_1\cdot \min\{\omega(t):t\in[a,b]\}$.  The convolution power $*^nf$ is supported on $[na,nb]$,
and  so  $\nm{*^nf}_1\le |*^nf|_1\cdot\max\{\omega(t):t\in[na,nb]\}$. We deduce that
$$
\nm{*^nf}_1\le {\nm f}_1^n\cdot\frac{\max\{\omega(t):t\in[na,nb]\}}{(\min\{\omega(t):t\in[a,b]\})^n}.
$$
Taking $n$th roots and using the spectral radius formula and the fact that
$\omega(t)^{1/t}\to 0$, we see
that $f$ is a quasinilpotent element of $L^1(\bR^+,\omega)$.
Since  the operator norm in $\cA$ is bounded by the $L^1$ norm $\Nm_1$, it follows that
 $M_f$ is
quasinilpotent in $\cA$ too.

The fact that $\omega(t)\to 1$ as $t\to 0$ ensures that for small $\veps>0$, the $L^1$ norm of any nonnegative function $f$ whose integral is 1 and which is supported on $[0,\veps]$, is close to 1. The corresponding operators $M_f$ form a contractive approximate identity for $\cA$. Also, $\cA$ has, for each $\alp\ge 0$, a ``standard ideal'' $J_\alp$ consisting of the norm closure of operators $M_f$ with $f\in L^1$ supported on $[\alp,\infty)$. We shall show:

\bthm\l cdvmain!
For any radical weight $\omega$, the algebra $\cA(\omega)$ is an integral domain with cai.
 If the radical weight $\omega$ satisfies Domar's criterion, then the standard ideals are
the only nonzero closed ideals of $\cA$.
\ethm

Note that $J_\alp\cdot J_\beta \subset J_{\alp+\beta}$, so
$J_\alp\ne \overline{J_\alp^2}$ for $\alp>0$.  In particular, the nontrivial standard ideals do not have any approximate identity and are not $r$-ideals. So when the above theorem is proved, we will have shown that the algebra $\cA$ has all the properties claimed at the start of the section, provided the weight $\omega$ satisfies Domar's criterion.

Recall that $L^1_\loc(\bR^+)$ denotes the Fr\'echet space of locally integrable measurable  functions on $\bR^+$.
For a function $f\in L^1_\loc(\bR^+)$,
 we define  $\alp(f)$ to be the minimum of the support of $f$ (or $+\infty$ if $f=0$).
We will use the  Titchmarsh convolution theorem (see e.g.\ \cite[Theorem 4.7.22]{Dal}),
which states for example that
$\alp(f*g)=\alp(f)+\alp(g)$ for  $f,g\in L^1_{\loc}(\bR^+)$.
 In particular this is the case when $f\in L^p(\bR^+,\omega)$ and $g\in  L^q(\bR^+,\omega)$ for some $p,q\ge 1$ (for on compact intervals
$\omega$ is bounded away from $0$ and so by the H\"older inequality it is clear that
$L^p(\bR^+,\omega)\subset L^1_\loc$ for every $p\ge 1$). For an operator $T\in\cA$ we define $\alp(T)=\inf\{\alp(Tf):f\in L^2(\bR^+,\omega)\}$.  

\blem\l alpha!
For each radical weight $\omega$ and each $S,T\in\cA(\omega)$, we have $\alp(ST)=\alp(S)+\alp(T)$.
\elem
\proof
Let $f\in L^2(\bR^+,\omega)$ and let $g=T(f)$, so that
 $\alp(g)\ge \alp(T)$. Write $g_1(x)=g(x)\cdot {\bf 1}_{x\le \alp(S)+\alp(T)}$.
Then $g_1$ is compactly supported with $|g_1|^2\in L^1$, and so  $g_1=\del_{\alp(T)}*g_0$
for some  $g_0\in L^2$. Then $S(g_1)=\del_{\alp(T)}*S(g_0)$,
 because $S$ is a norm limit of convolution operators.
Hence  $\alp(S(g_1))=\alp(T)+\alp(S(g_0))\ge \alp(S)+\alp(T)$.
The functions $Sg_1(x)$ and $Sg(x)$ agree for $x\le \alp(T)+\alp(S)$,
 because $g(x)$ and $g_1(x)$ agree for such values, and the subspace
of functions $f$ with $\alp(f)\ge \alp$  is
invariant for all operators under consideration.
Thus $\alp(STf)=\alp(Sg)\ge \alp(T)+\alp(S)$, so for all $S$ and $T$ we have $\alp(ST)\ge\alp(S)+\alp(T)$.

To prove the converse, we will use the fact that the function $f \mapsto \alp(f)$ is
upper semicontinuous on $L^p(\bR^+,\omega)$.  This is  because if the minimum of the
support of $f\in L^p$ is $\alp$, then for $\veps>0$ the integral
$\int_\alp^{\alp+\veps} \, |f(t)|^p \, dt$ is strictly positive.
Hence for functions $g$ sufficiently close to $f$ in $p$ norm, we will have
$\int_\alp^{\alp+\veps} \, |g(t)|^p \, dt>0$ also, and in particular $\alp(g)<\alp+\veps$.
Since the set $C_{00}(\bR^+)$ of continuous functions
of compact support is dense in $L^2$,
given $\veps>0$ we may pick $f,g\in C_{00}(\bR^+)$ such that $\alp(Sf)\le\alp(S)+\veps/2$ and
 $\alp(Tg)\le\alp(T)+\veps/2$. Then $ST(f*g)=S(f*(Tg))$,
 because $T$ is a norm limit of convolution operators.  This is equal to $S(f)*T(g)$
because $S$ is a norm limit of convolution operators. By the Titchmarsh
convolution theorem, $\alp(S(f)*T(g))=\alp(Sf)+\alp(Tg)$, and  so
$$\alp(ST)\le \alp(ST(f*g))=\alp(Sf)+\alp(Tg)\le \alp(S)+\alp(T)+\veps.$$
Hence $\alp(ST)\le\alp(S)+\alp(T)$, and the lemma is proved.\endproof

\bcor\l id!
The algebra $\cA(\omega)$ is an integral domain.
\ecor

We have defined the standard ideal $J_\alpha$ to be the closure in $\cA$ of the operators $M_f$ with $f\in L^1(\bR^+,\omega)$ and $\alp(f)\ge\alp$.
There is another obvious closed ideal,
namely  $I_\alpha=\{T\in\cA:\alp(T)\ge \alp\}$.  We now
 show that these two ideals coincide.

\blem\l same!
Let $\omega$ be any radical weight. Then for each $\alp\ge 0$, the ideals $I_\alp$ and $J_\alp\subset\cA(\omega)$
are the same.
\elem
\proof Every operator $M_f$ with $f\in L^1(\bR^+,\omega)$ and $\alp(f)\ge\alp$ is plainly in $I_\alp$, which is  closed because it is the set of operators in $\cA$ which map $H$ into the closed subspace of functions supported on $[\alp,\infty)$. Therefore $J_\alp\subset I_\alp$.  Conversely, let $T\in I_\alp$ with $T=\lim_iM_{f_i},\ f_i\in L^1(\bR^+,\omega)$. We claim that $T\in J_\alp$, which will imply
 that the two ideals are the same. To prove this, let $f$ be smooth and compactly
supported in $[0,\infty)$. The operator $T\cdot M_f=\lim_i M_{f_i*f}$, and the function $\gam=T(f)=\lim_i f_i*f$,
are supported on $[\alp,\infty)$.  Thus $(f_i*f)\cdot {\bf 1}_{[0,\alp]}\to 0$ in
$L^2([0,\alp],\omega)$, and even in  $L^1([0,\alp],\omega)$, because the $L^1$ norm on the compact set $[0,\alp]$ is bounded by a constant times the $L^2$ norm.
 It follows that if $\gam_i=(f_i*f)\cdot {\bf 1}_{[0,\alp]}\in\cA$ then the operator norm $\nm{\gam_i}\to 0$ (for the operator norm is at most the $L^1$ norm,
which is known to tend to zero). So $TM_f=\lim_iM_{(f_i*f)\cdot {\bf 1}_{(\alp,\infty)}}\in J_\alp$.
However,
 the algebra $\cA$ is known to have a sequential cai consisting of (convolution operators by) functions $u_i$ which are smooth and compactly supported in $[0,\infty)$. So $T=\lim_i TM_{u_i}\in J_\alp$ also, and $I_\alp=J_\alp$.\endproof

\blem\l x1!
If the radical weight $\omega$ satisfies Domar's criterion, then for each $t>0$ the
integral $\int_0^\infty \, (\omega(x+t)/\omega(x))^2 \, dx$ is finite.
\elem
\proof The function $\eta(x)=-\log\omega(x)$ is a convex function by Domar's criterion.
Thus the ratio $\omega(t+x)/\omega(x)=\exp(\eta(x)-\eta(t+x))$ is a decreasing function of $x$.
We know that $\eta(x)\ge x^{1+\veps}$ for large $x$, so that
$\eta(nt)\ge t^{1+\veps}n^{1+\veps}$ for large enough $n$. Since $\eta(0)=0$,
 we have $\summ r1n\eta(rt)-\eta((r-1)t)\ge t^{1+\veps}n^{1+\veps}$, and so the largest term $\eta(nt)-\eta((n-1)t)$ must dominate $n^\veps t^{1+\veps}$.
This implies that for all large enough $n$ we have
$\omega(nt)/\omega((n-1)t)\le e^{-cn^\veps}$, where $c=t^{1+\veps}$.
 Since $\omega(x+t)/\omega(x)$ is a decreasing function, this implies that
  $\int_0^\infty \, (\omega(x+t)/\omega(x))^2 \, dx < \infty$. \endproof

\bcor\l x2!
If $\omega$ satisfies Domar's criterion, $f\in L^2(\bR^+,\omega)$ and $t>0$ then $f*\del_t\in L^1(\bR^+,\omega)$.
\ecor

\begin{proof} By the Cauchy-Schwartz inequality we have
$$\int_0^\infty \, |f(x)| \, \omega(x+t) \, dx \le (\int_0^\infty
\, |f(x)|^2 \, \omega(x)^2 \, dx)^{1/2} \; (\int_0^\infty \,
(\omega(x+t)/\omega(x))^2 \, dx)^{1/2}<\infty .$$ \end{proof}

\blem\l x3!
Suppose a radical weight $\omega$ satisfies Domar's criterion, $T\in\cA(\omega)$ is nonzero, and $\alp(T)=\alp$. Then the closed principal ideal $\overline{T\cdot \cA}$ is equal to $I_\alp$.
\elem
\proof By \reflem{alpha} we have $\overline{T\cdot \cA}\subset I_\alp$. Conversely, since $I_\alp=J_\alp$ it is enough to show that for each $f\in L^1(\bR^+,\omega)$ with $\alp(f)\ge \alp$ we have $M_f\in \overline{T\cdot \cA}$.
Any such $f$ is a norm limit of functions  $f_n \in L^1(\omega)$ with $\alp(f_n)>\alp$;
and so it is enough to show that
 $M_f\in \overline{T\cdot \cA}$ when $\alp(f)>\alp$. Given such a function
 $f$, pick $g\in L^2(\bR^+,\omega)$ with $\alp(Tg)<\alp(f)$ (this is
possible because the infimum of values $\alp(Tg)$ is by hypothesis equal to $\alp$),
 and pick  $t>0$ such that we have $t+\alp(Tg)=\alp(\del_t*Tg)<\alp(f)$. By \refcor{x2}\  the function $h_0=\del_t*g$ is
in  $L^1(\bR^+,\omega)$, as also is $h=\del_t*Tg$.

By Domar's theorem the closed ideal generated by $h$ in $L^1(\bR^+,\omega)$ is standard. It therefore contains every
function $k$ in $L^1(\bR^+,\omega)$ with $\alp(k)\ge \alp(h)$, so it contains the function $f$. There is a sequence of
 functions $u_i\in L^1(\bR^+,\omega)$ with $u_i*h\to f$. The operator norm on $\cA$ is bounded by the $L^1$ norm so $M_{u_i*h}\to M_f$ in $\cA$.
Now $T$ is a norm limit of convolution operators:  $T=\lim_j \, M_{\tau_j}$ with
$\tau_j\in L^1(\bR^+,\omega)$.  So for any $\gam\in C_{00}(\bR^+)$, we have
\begin{align*}
M_{u_i*h}(\gam) &= u_i*(\del_t*Tg)*\gam=u_i*\del_t*(\lim_j\tau_j*g)*\gam \\
 &=\lim_j\tau_j*u_i*\del_t*g*\gam=T\cdot M_{u_i*h_0}(\gam),
\end{align*}
because convergence of $\tau_j*g$ occurs in $L^2(\bR^+,\omega)$, and both $u_i*\del_t$ and $\gam$ are in $L^1(\bR^+,\omega)$, and $L^1(\bR^+,\omega)$ acts continuously on $L^2(\bR^+,\omega)$
by convolution.

Since $C_{00}$ is dense in $L^2$,
the operators $T\cdot M_{u_i*h_0}$ and $M_{u_i*h}$ are equal.
Hence the operator  $T\cdot  M_{u_i*h}$ is in the principal ideal $T\cA$. Therefore the
closure $\overline{T\cA}$ contains $M_f$ for every $f\in L^1(\bR^+,\omega)$ with $\alp(f)>\alp$. Thus  $\overline{T\cA}=I_\alp$ as claimed.  \endproof

\smallskip

\ni {\em Proof of \refthm{cdvmain}:} By \refcor{id}, $\cA(\omega)$
is an integral domain for any radical weight $\omega$. Let
$J\subset\cA$ be any nonzero closed ideal and let
$\alp=\alp(J)=\inf\{\alp(T):T\in J\}$. We claim that $J=I_\alp$.
From the definition it is plain that $J\subset I_\alp$. For the
converse, choose $T_n\in J$ with $\alp(T_n)\to\alp$.  By
\reflem{x3}, $J$ contains $I_\beta$ for a sequence of values $\beta$
tending to $\alp$. In particular $J$ contains the operator $M_f$ for
every $f\in   L^1(\bR^+,\omega)$ with $\alp(f)>\alp$. The closure of
this set includes every $M_g$ with $\alp(g)\ge \alp$. Hence  it
contains $J_\alp=I_\alp$. Thus $J=I_\alp$ as claimed. $\Box$

\smallskip

It is easy to see that the algebras $A(\omega)$ above contain
no idempotents, and are `modular annihilator algebras'.  As in
\cite{Dal}, the function $u = 1$ in $L^1(\bR^+,\omega)$ corresponds to a single generator
for the algebra $A$.

\section{Pre-images of HSA's}

If $J$ is  a closed ideal in an approximately unital
operator algebra, we examine the relation between
${\mathfrak F}_A$ and ${\mathfrak F}_{A/J}$.  From Meyer's theorem (\cite{Mey},
\cite[Theorem 2.1.13]{BLM}) one can see that the `image' of ${\mathfrak F}_A$
in $A/J$ is a subset of ${\mathfrak F}_{A/J}$.

\begin{proposition} \label{ith}  If $J$ is a closed ideal in an
operator algebra $A$, and if $J$ has a cai, then
$q({\mathfrak F}_A) = {\mathfrak F}_{A/J}$, where $q : A \to A/J$ is the
canonical map.
\end{proposition} \begin{proof}  Indeed
suppose that $x \in A/J$ with
  $\Vert 1 - x \Vert \leq 1$ in  $A^1/J \cong (A/J)^1$.
 Since $J$ is an $M$-ideal  in $A^1$ (see e.g.\  \cite[Theorem 4.8.5]{BLM}),
it is proximinal \cite{HWW}.  Hence there is an element $z = \lambda 1 + a$ in
Ball$(A^1)$, with $\lambda \in \Cdb, a \in A$, such that
$\lambda 1 + a +J = 1 - x$.  It is easy to see now that $\lambda = 1$,
and $a+J = -x$.  Let $y = -a$.
 Then $\Vert 1 - y \Vert = \Vert 1 + a \Vert = \Vert z \Vert \leq 1$,
   so $y \in {\mathfrak F}_A$, and $q(y) = x$.
\end{proof}

\begin{proposition} \label{hio}
If $J$ is a closed ideal in an 
operator algebra $A$, and if $J$ has a cai too, then any closed
approximately unital subalgebra
 $D$ in $A/J$ is
the image of a closed approximately unital subalgebra
in $A$, under
the quotient map $q_J$ from $A$ onto $A/J$.
In fact $q_J^{-1}(D)$ will serve here.
\end{proposition}

\begin{proof}   The idea for this proof was found independently
by M. Almus.   Note that $J$ is an approximately unital  ideal
in $q_J^{-1}(D)$.  Moreover, $q_J^{-1}(D)/J \cong D$, which is
approximately unital.  So $q_J^{-1}(D)$ is
approximately unital by \cite[Proposition 3.1]{BR}.
 Another proof
follows immediately from 3.4 in
\cite{BR}, since,  in the language there,
$B \oplus_C C' = \beta^{-1}(C')$ clearly.
\end{proof}

\begin{corollary} \label{a4}  Let $J$ be a closed approximately unital
two-sided ideal in an operator algebra
$A$, and let $q_J : A \to A/J$ be the quotient map. \begin{enumerate}
\item [{\rm (i)}]  The open projections
in $(A/J)^{**}$ are exactly the $q_J^{**}(p)$, for
open projections $p$ in $A^{**}$.
\item [{\rm (ii)}] The HSA's in $A/J$ are precisely the images
of the HSA's in $A$, under  $q_J$.
 \item [{\rm (iii)}] The r-ideals in $A/J$ are
precisely the images of the r-ideals in $A$, under $q_J$.
\item [{\rm (iv)}]   An r-ideal (resp.\ HSA) in $A/J$ of the
form $\overline{x (A/J)}$ (resp.\ $\overline{x (A/J) x}$) for some $x \in {\mathfrak F}_{A/J}$,
is the image of an r-ideal (resp.\ HSA) in $A$ of the
form $\overline{y A}$ (resp.\ $\overline{y A y}$) for some $y \in {\mathfrak F}_{A}$.
 \end{enumerate}
 \end{corollary}

\begin{proof}   It is easy to see that  the images of
HSA's (resp.\ r-ideals) in $A$,
are HSA's (resp.\ r-ideals) in $A/J$.   If $p$ is open in $A^{**}$
then $p$ is the weak* limit of a net $(a_t)$ in $A$ with $a_t = p a_t p$.
Then $q_J^{**}(p)$ is the weak* limit of a similar net in $A/J$, so is open
there.   Items (i)--(iii) follow easily from these observations,
and Proposition \ref{hio}.  For (i), if $p$ is the support projection of
 $D' = q_J^{-1}(D)$, where $D$ is the HSA associated with $p$, then $q_J^{**}(p)$ is the support projection of $D$.  So the
 open projections in $(A/J)^{**}$ are precisely the $q_J^{**}(p)$,
 for open projections $p \in A^{**}$.

   Item  (iv) follows easily from Proposition \ref{ith};
and that result also leads to another
proof of (i)-(iii), which seems to give a possibly
different preimage.   We give the argument in the r-ideal case:
  Let $K$ be an r-ideal in $A/J$.
By Corollary \ref{rd1}, there is
a lcai $(e_t)$ in $K$ with $\Vert 1 - 2 e_t \Vert \leq 1$.  As above
we obtain
$x_t \in {\mathfrak F}_A$, with $q_J(x_t) = e_t$.
The closure of the sum
   of the right ideals $x_t A$, is an r-ideal $K'$ in $A$
   by     Theorem \ref{corea11}.  Moreover $q_J(K')$
   is contained in the closure of the union
   of the  $q_J(x_t) (A/J) = e_t (A/J) \subset K$.
 Conversely, for any $a \in A$, we have $e_t (a+J) =  q_J(x_t a)$;
 and since  $(e_t)$ is a lcai for $K$ it follows that $K \subset
 q_J(K')$.  So $q_J(K') = K$.
\end{proof}

This technique seems applicable to  other `constructions' besides quotients, such as direct limits, ultrapowers, interpolated
operator algebras, etc.  See \cite[Sections 2.2 and 2.4]{BLM} for some of
 these constructions.  Indeed the results apply directly to ultraproducts
 because they are quotients of the type described in this section.

\bigskip

\section{Other constructions of $r$-ideals}

The following is an improvement of \cite[Proposition 3.1]{ABS}, and
also answers the question in the Remark following it.

\begin{theorem}  \label{best}  If $A$ is an operator algebra
with left cai,
which is a left ideal in an operator algebra $B$, then $\overline{(1-x) A}$
is an r-ideal in $A$ for all $x \in {\rm Ball}(B)$.  \end{theorem}

 \begin{proof}  We may assume that $B$ is unital.
 Certainly $J = \overline{(1-x) A}$ is a right ideal, the question
is whether it has a left cai.  Note that
$e_n = 1 - \frac{1}{n} \sum_{k=1}^n \, x^k$ defines a bounded net,
$\Vert 1 - e_n \Vert \leq 1$, and $e_n (1-x) = 1-x - \frac{1}{n} (1 - x^{n+1})
\to 1-x$.  Suppose that $(f_t)$ is a left cai for $A$, with weak* limit $f
 \in A^{\perp \perp}$, which is a projection and a left identity
  for $A^{\perp \perp}$.  We may view  $(e_n f_t)$ as a net in $(1-x)  A \subset J$, with the product indexing,
and it is easy to see using the above that $e_n f_t (1-x) a \to (1-x) a$ for all $a \in A$.   Suppose that a subnet $((1 - e_{n_\mu}) f_{t_\mu})$ converges
weak* to an element $r$.  Then it is easy to see that $r$ is a contraction
in $A^{\perp \perp}$, so that $fr = r$.  Also, $e_{n_\mu} f_{t_\mu} \to f-r$, so that $f-r \in J^{\perp \perp}$.  Since $e_{n_\mu} f_{t_\mu} (1-x) a \to (1-x) a$, we
have $(f-r) z = z$ for all $z \in J$, hence for all $z \in J^{\perp \perp}$.  So $f-r$ is a left identity for $J^{\perp \perp}$, hence
it is idempotent.  That is, $f - fr - rf + r^2 = f-r$, which by
a fact above implies that $r^2 = rf$.  If we choose $(f_t)$
 so that $f_s f_t \to f_t$ with $t$ (as in Corollary \ref{rd1} above,
 or \cite[Corollary 2.6]{BHN}),
 we may assume that
$f_t f = f_t$, which forces $r f = r$ by definition of $r$.
Thus $r$ is idempotent, hence is a projection, and so $f-r$ is
a projection too.  By e.g.\ \cite[Proposition 2.5.8]{BLM},
$J$ has a left cai.
 \end{proof}

{\bf Remarks.}  1) \ We do not have a clean formula for the left cai in the last result,
although there is one for a left bai: the net $(e_n f_t)$ in the proof
is a left bai.
This illustrates the fact that although we may know that left  cai
exist of a nice form (as in Theorem \ref{read}
or in the Remark after Theorem \ref{sephs}), we may not be able to
write a simple expression for them.

\smallskip

2) \ Considering the example of the $2\times 2$ matrices supported on the first
column, shows that the last result is best possible.  That is,
the hypothesis of a left cai is not removable.

\bigskip

The following result is actually equivalent to the last theorem:

\begin{corollary}  \label{abm}  If $A$ is an operator algebra
with left cai,
and if $\eta : A \to A$ is a completely contractive
 left $A$-module map, or if
$\eta \in {\rm Ball}(A^{**})$ satisfies $\eta A \subset A$,
then $\overline{(1-\eta) A}$ is an r-ideal of $A$.
 \end{corollary}

\begin{proof}
Consider $B = \{ \eta \in A^{**} : \eta A \subset A \}$,
an operator algebra containing $A$ as a left ideal.   Thus the second case of
our result follows from Theorem  \ref{best}.

The set of completely bounded
 left $A$-module maps `equals' $B$ by \cite[Theorem 6.1]{Bonesided}
 (note that a hypothesis in the latter theorem was removed
 in  \cite[Corollary 2.6]{BHN}), and hence this case follows
 by the last paragraph.
 \end{proof}

{\bf Remarks.}  1) \ By another equivalence in \cite[Theorem 6.1]{Bonesided},
the last result is correct with $\eta$ a contraction
in the operator space
left multiplier algebra ${\mathcal M}_{\ell}(A)$ (see \cite[Chapter 4]{BLM}
for the definition of the latter).  If $A$ is approximately unital
then ${\mathcal M}_{\ell}(A) = LM(A)$, the ordinary left multiplier algebra.

\smallskip

2)  The first result of this type that we are aware of dates to
 2005 (see \cite[Lemma 6.8]{BHN}, but this is much less general).  See
 also \cite{KLU}, for some recent Banach algebra variants.

\bigskip

According to \cite[Corollary 2.7]{BHN}, there is a
bijective correspondence
between the classes of r-ideals, $\ell$-ideals, and HSA's,  of $A$.
   One may ask what is the
 $\ell$-ideal and HSA matching the r-ideal in Theorem  \ref{best},
 in terms of $x$?   In general we do not have a simple answer.
 However we have:

\begin{proposition} \label{isb}  If $A$ is an operator algebra
with cai,
which is an ideal in an operator algebra $B$, and $x \in {\rm Ball}(B)$, then the
$\ell$-ideal and HSA matching the r-ideal $\overline{(1-x) A}$,
are $\overline{A(1-x)}$ and $\overline{(1-x) A (1-x)}$.
\end{proposition}

\begin{proof}  We may assume that $B$ is unital.  Then $A$ corresponds to
a central projection $p$ in $B^{**}$, whereas the r-ideal $J = \overline{(1-x)B}$
in $B$ has a support projection $e \in B^{**}$, say.  Then
$J^{\perp \perp} = e B^{**},  \overline{B(1-x)}^{\perp \perp} = B^{**} e$,
and $\overline{(1-x)B(1-x)}^{\perp \perp} = e B^{**} e$ by facts in
\cite{BHN}.  Since
$e$ and $p$ commute, we have that
$e B^{**}  \cap B^{**} p = ep B^{**}$.  By \cite[5.2.7]{BZ},
$A^\perp + J^\perp$ is closed, and by a formula in the proof of
\cite[5.2.9]{BZ}, we have
that $$(A \cap J)^{\perp \perp} = (A^\perp + J^\perp)^\perp
= A^{\perp \perp} \cap J^{\perp \perp} =
 p B^{**}  \cap e B^{**} = ep B^{**}. $$
Similarly, $(A \cap \overline{B(1-x)})^{\perp \perp} = B^{**} ep$.
Now $\overline{A(1-x)} \subset A \cap \overline{B(1-x)}$.
 Conversely if $z \in A \cap \overline{B(1-x)}$ then since
$A$ has a cai $(e_t)$ and $A$ is an ideal in $B$, we have
$z = \lim_t \, e_t z \in \overline{A(1-x)}$.
Thus $A \cap \overline{B (1-x)} = \overline{A(1-x)}$,
and, similarly, $A \cap \overline{(1-x) B} = \overline{(1-x) A}$.
It follows that $\overline{A(1-x)}$ is the $\ell$-ideal matching $\overline{(1-x) A}$.
By \cite[Corollary 2.8]{BHN}, the corresponding HSA will be the intersection of these,
which also equals
their product, which can be seen to be
$\overline{(1-x) A (1-x)}$, using the cai for $A$.  We do not need this
here, but this HSA also equals
$A \cap \overline{(1-x)B(1-x)}$, since the latter equals $A \cap \overline{(1-x)B}  \cap \overline{B(1-x)}
= \overline{(1-x)A} \cap \overline{A(1-x)} . $
\end{proof}

{\bf Remark.}
 The result is not true if $A$ only has a one-sided cai.  For example
if $x = E_{11},$ and $A = C_2$ as in the example after Theorem \ref{best}.

\section{Positive maps between operator algebras}

The {\em size} of ${\mathfrak F}_A$, and what all it contains,
 seems to be an important and possibly quite difficult
question for nonunital operator algebras $A$. Of course for  unital
$A$ the answer is trivial.  Note too that for a $C^*$-algebra $A$
with positive cai $(e_t)$, then one obtains a probably quite good
idea of what is contained in ${\mathfrak F}_A$, by meditating on the
simple fact that $a^2 + a x a \in {\mathfrak F}_A$ for all $x \in
{\rm Ball}(A)$ and $a \in {\rm Ball}(A)_+$ (this follows
since the product $y \; {\rm diag}(1, x) \; y^*$ is a contraction
where $y = [\sqrt{1-a^2} \; \; \; a ]$).  In particular,
 $e_t^2 + e_t x e_t \in {\mathfrak F}_A$ in this case.
Unfortunately, this seems to be far from true
for nonselfadjoint operator algebras.

\medskip

{\bf Remark.}  We remark in passing that the only idempotents that could be contained in
$\Rdb^+ {\mathfrak F}_A$, are orthogonal projections.  Also, note that ${\mathfrak F}_A$ can contain selected unitaries (eg.\ certain
functions valued in the unit circle on certain subsets of $[0,1]$), but not nonunitary
isometries (by e.g.\ Corollary \ref{perm2}).

\medskip

\begin{lemma} \label{refe}  If $A$ is an approximately unital
operator algebra, then ${\mathfrak F}_A$ is  weak* dense
in ${\mathfrak F}_{A^{**}}$.
\end{lemma}

\begin{proof} (We are indebted to the referee for supplying this proof.)
Assume that $A \subset B(H)$.
Let $(v_t)$ be a cai as in Theorem \ref{newr},
with numerical range in the wedge
shape region of angle $2 \rho$ described there,
where $\rho  \to 0$  as $t$ increases through the directed set.
 Write $v_t = a_t + ib_t$,
for selfadjoint $a = a_t$ and $b = b_t$.
Because of the position of the numerical range of $v_t$,
 $a$ is a positive contraction.
Also, for all states
$\varphi$ on a $C^*$-algebra generated by $A$,
 we have $|\varphi(b)| \leq  (\tan \rho) \varphi(a)$. So
$a \tan \rho \pm b \geq 0$.  By a well known fact about selfadjoint operators
(which is a
pleasant exercise to prove), there exists a selfadjoint $c \in B(H)$
with $b = a^{\frac{1}{2}} c a^{\frac{1}{2}}$ and $\Vert c \Vert \leq \tan \rho$.
   Then $v_t = a^{\frac{1}{2}} (1 + ic) a^{\frac{1}{2}}$.
Setting $c_t = c$ we have $c_t \to 0$ with $t$.

Let $z \in {\rm Ball}(A^{**})$; by Goldstine's lemma we may
choose $z_i \in {\rm Ball}(A)$ with $z_i \to z$ weak*.
Fix $\delta > 0$ and set $x_{i,t,\delta} = (1-\delta) v_t + (1-2\delta)
v_t z_i v_t$.  It is easy to check that
$1-x_{i,t,\delta} = 1 - a + a^{\frac{1}{2}} w
 a^{\frac{1}{2}},$
where $$w = \delta 1 - (1-\delta) ic
- (1-2\delta) (1 + ic) a^{\frac{1}{2}} z_i a^{\frac{1}{2}} (1 + ic).$$
Since $c = c_t \to 0$ with $t$, for $t$ `large' we have
$$\Vert w
\Vert \leq \delta + (1-\delta) \Vert c \Vert + (1-2\delta) \Vert 1 + ic \Vert^2 \leq 1 .$$
Hence $\Vert 1 - x_{i,t,\delta} \Vert \leq 1$,
 since the product $y \; {\rm diag}(1, w) \; y^*$ is a contraction
where $y = [\sqrt{1-a} \; \; \; a^{\frac{1}{2}} ]$.
Thus $x_{i,t,\delta} \in {\mathfrak F}_A$.  Since $v_t z_i v_t \to z_i$ in norm with
$t$, it follows that
$(1-\delta) 1 + (1-2\delta) z_i$ is in the weak* closure of
${\mathfrak F}_A$ for every $i$.  Hence $1 + z$ is in this weak* closure
too.  \end{proof}

Below we will also consider {\em unital operator spaces}: subspaces
$A$ of $B(H)$ containing $I_H$ (see \cite{BN} for a matrix norm
characterization of these).   Here ${\mathfrak F}_A = \{ x \in A :
\Vert 1_A - x \Vert \leq 1 \}$.  One may define a cone in any
operator algebra (or unital operator space) $A$ by considering
${\mathfrak c} = {\mathfrak c}_A = \Rdb^+ {\mathfrak F}_A$. Probably
$\frac{1}{2} {\mathfrak F}_A$ should be considered to be the
analogue of the positive part of the {\em unit ball} of a
$C^*$-algebra. Similarly, one obtains cones ${\mathfrak c}_n$ in
$M_n(A)$ for every $n \in \Ndb$.

The following shows that ${\mathfrak c}_A$ is large enough to determine $A$:

\begin{corollary} \label{hi}  Suppose that $A$ and $B$ are approximately
unital closed subalgebras of $B(H)$, or unital subspaces of $B(H)$
with identities $1_A$ and $1_B$
corresponding to projections on $H$.  If
${\mathfrak c}_A \subset {\mathfrak c}_B$ then $A \subset B$.  Hence
$A = B$ iff ${\mathfrak c}_A = {\mathfrak c}_B$.
 \end{corollary}  \begin{proof}   First assume that
$A$ and $B$ are unital.  If $x \in {\rm Ball}(A)$ then
$1_A$ and $1_A + x$ are in ${\mathfrak F}_A \subset
{\mathfrak c}_B$, and so $1_A, x \in B$.  Hence $A \subset B$.

In the general case, taking weak*-closures in
$B(H)^{**}$, we have by
Lemma \ref{refe}  that ${\mathfrak F}_{A^{\perp \perp}} =
\overline{{\mathfrak F}_A}^{w*} \subset \overline{{\mathfrak F}_B}^{w*} =
{\mathfrak F}_{B^{\perp \perp}}$. By the last paragraph,
$A^{\perp \perp} \subset B^{\perp \perp}$, and hence
$A = A^{\perp \perp} \cap B(H) \subset B = B^{\perp \perp} \cap B(H)$. \end{proof}

\begin{definition} \label{defcp}
We say that a map $T : A \to B$ between  operator algebras, or
between unital operator spaces, is {\em operator completely
positive}, or {\em OCP}, if there is a constant $C > 0$ such that
$T_n({\mathfrak F}_{M_n(A)}) \subset C {\mathfrak F}_{M_n(B)}$ for
every $n \in \Ndb$. We study these maps below.  If $A$ and $B$ are
operator algebras, but not unital, then we will also require $T$ to
be completely bounded (this is automatic if $A$ is unital).
\end{definition}

Some remarks on Definition \ref{defcp}: First,  the definition  is
`positive homogeneous' in $C$.  That is, $T : A \to B$ satisfies
$T_n({\mathfrak F}_{M_n(A)}) \subset C {\mathfrak F}_{M_n(B)}$, iff
$R_n({\mathfrak F}_{M_n(A)}) \subset {\mathfrak F}_{M_n(B)}$ where
$R = \frac{T}{C}$.  Thus we may usually assume that $C = 1$. Second,
we will also use the fact that $x \in {\mathfrak c}$ iff there is a
constant $C > 0$ with $x + x^* \geq C x^* x$.  Third, it is obvious
that a completely contractive unital linear map between unital
operator spaces is OCP. Finally, we remark that if $\varphi : A \to
B$ is a completely contractive homomorphism between operator
algebras, then $\varphi$ is OCP. Indeed by Meyer's theorem
(\cite{Mey}, \cite[Theorem 2.1.13]{BLM}), we can extend $\varphi$ to
a completely contractive unital homomorphism between unitizations,
and then the result is obvious by the third remark.

We write $C^*(A)$ for a $C^*$-algebra that contains $A$ completely
isometrically as a subalgebra if $A$ is an operator algebra, or as a
unital subspace if $A$ is a unital operator space (with $1_A =
1_{C^*(A)}$ in this case), and which is generated by $A$.

\begin{lemma}   \label{ik}  If $\varphi : A \to B(H)$ is a map from an
operator algebra, or from a unital operator space, that extends to a
completely positive map from $C^*(A)$ into $B(H)$, then $\varphi$ is
OCP.
\end{lemma}

\begin{proof}
 We may assume without loss of generality that $A$ is a $C^*$-algebra.
Suppose that $x \in {\mathfrak F}_{A}$. Then $$\varphi(x) +
\varphi(x)^* = \varphi(x + x^*) \geq \varphi(x^* x) \geq C
\varphi(x)^* \varphi(x) ,$$ for a constant $C = \Vert \varphi
\Vert_{cb}^{-1}  > 0$, by the
Kadison-Schwarz inequality (see e.g.\ \cite{SOC,Pn}). Thus $\varphi({\mathfrak F}_A)
\subset \Vert \varphi
\Vert_{cb} \,  {\mathfrak F}_B$.  Similarly for matrices.  So  $\varphi$ is
OCP.
\end{proof}

\begin{lemma}   \label{ikhuh}  If $A$ is a $C^*$-algebra or operator system,
then $x \in {\rm Ball}(A)_+$ iff $zx \in {\mathfrak F}_{A}$ for
all $z \in {\mathfrak F}_{\Cdb}$.
\end{lemma}

\begin{proof}   ($\Rightarrow$) \ Left to the reader.

($\Leftarrow$) If $x$ satisfies this property then for any $z \in {\mathfrak F}_{\Cdb}$, we have
$|z|^2 x x^* \leq 2 \, {\rm Re}(zx)$, so that
${\rm Re}(z \langle x \zeta, \zeta \rangle) \geq 0$ for any unit vector
$\zeta \in H$.
It is a pleasant exercise in calculus that if the latter holds for
all $z \in {\mathfrak F}_{\Cdb}$ then $\langle x \zeta, \zeta \rangle
\geq 0$.  So $x$ is positive, and it is easy to see that it has to be
a contraction.
\end{proof}

\begin{theorem}   \label{ik3}  If $T : A \to B$ is a map between
$C^*$-algebras or operator systems
then $T$ is completely positive iff $T$ is OCP.
\end{theorem}

\begin{proof}
By virtue of Lemma  \ref{ik}  we need only prove one direction.
Suppose that $T$ is OCP.
By one of the observations
below Definition \ref{defcp},  we may assume that  $C = 1$ in the definition of OCP.
  If $x \in {\rm Ball}(A)_+$ and $z \in {\mathfrak F}_{\Cdb}$
then $zx  \in {\mathfrak F}_{A}$
by Lemma \ref{ikhuh}.  Thus $zT(x) = T(zx)  \in {\mathfrak F}_{B}$,
and so
$T(x) \geq 0$ by Lemma \ref{ikhuh}.  A similar argument applies to matrices.
\end{proof}

\begin{theorem}   \label{kk}  If $T : A \to B(H)$ is an OCP
map on a  unital operator space $A$, then  the canonical extension
$\tilde{T} : A + A^* \to B(H): x + y^* \mapsto T(x) + T(y)^*$ is
well-defined and completely positive.
\end{theorem}

\begin{proof}  As in the last proof,
we  may assume that  $C = 1$ in the definition of OCP. In this case
notice that by the last result applied to the restriction
of $T$ to $\Cdb 1$, we have $0 \leq T(1) \leq I$.
Assume first that $\varphi : A \to \Cdb$ is OCP. Since $|1
-\varphi(1) - \varphi(x)| \leq 1$ for all $x \in {\rm Ball}(A)$, we
have $1 -\varphi(1) + |\varphi(x)| \leq 1$, so that $\Vert \varphi
\Vert \leq \varphi(1)$. Hence $\Vert \varphi \Vert = \varphi(1)$.
Thus $\varphi$ extends by the Hahn-Banach theorem to a functional
$\psi : A + A^* \to \Cdb$ satisfying $\Vert \psi \Vert = \psi(1)$.
The latter implies that $\psi$ is positive \cite{SOC,Pn}.

To see that  $\tilde{T}$ is well-defined,
 notice that if $x + y^* = 0$, and if $\varphi$ is any state on $B = B(H)$,
 then by the last paragraph,  $\varphi \circ T$ extends to a positive map
 $\psi$ on $A + A^*$, so that $\varphi(T(x) + T(y)^*) = \psi(x+y^*) =
 0$.  Since this holds for every state on $B$ we have $T(x) + T(y)^* = 0$.

Similarly, if $x+y^* \geq 0$ then $\varphi(T(x) + T(y)^*) = \psi(x+y^*)
 \geq 0$.  Since this holds for every state on $B$ we have $T(x) + T(y)^*
 \geq 0$.  Thus $\tilde{T}$ is positive on $A + A^*$.  Applying
 this at every matrix level to $T_n$, we see that $\tilde{T}$ is
 completely positive on $A + A^*$.  \end{proof}

\begin{lemma}   \label{ik2}   If $T : A \to B(H)$ is an OCP map on
an approximately unital operator algebra, and if ${\mathfrak F}_{M_n(A)}$
is weak* dense in ${\mathfrak F}_{M_n(A^{**})}$ for
all $n \in \Ndb$, then the canonical weak*
continuous extension $\tilde{T} : A^{**} \to B(H)$ on the unital
operator algebra $A^{**}$ is OCP.
\end{lemma}

\begin{proof}   As in the last proofs,
  we may assume that  $C = 1$ in the definition of OCP.
Suppose that $\eta \in {\rm Ball}(A^{**})$.  By hypothesis,
there exists $(y_\lambda) \subset
 {\mathfrak F}_A$, with $y_t \to 1 + \eta$ weak*.  Then  $\Vert 1 - T(y_t) \Vert
\leq 1$, and in the weak* limit, $\Vert 1 - \tilde{T}(1 + \eta) \Vert
\leq 1$.  A similar argument prevails at the matrix level, so that
$\tilde{T}$ is OCP.
\end{proof}

\begin{theorem}   \label{ikc}  {\rm (Extension and Stinespring
dilation for OCP maps)} \ If $T : A \to B(H)$ is a
map on a unital operator space or on an approximately unital
operator algebra, and if $B$ is a $C^*$-algebra
containing $A$, then $T$ is OCP iff $T$ has a completely positive
extension $\tilde{T} : B \to B(H)$. This is equivalent to being able
to write $T$ as the restriction to $A$ of $V^* \pi(\cdot) V$ for a
$*$-representation $\pi : B \to B(K)$, and an operator $V : H \to
K$.  Moreover this can be
done with $\Vert T \Vert = \Vert T \Vert_{\rm cb} = \Vert V \Vert^2$, and
this equals $\Vert T(1) \Vert$ if $A$ is unital.
\end{theorem} \begin{proof}  As before, we only need prove one
direction of the first `iff'.  If $T$ is OCP and $A$ is a unital operator space,
then by the last result we can extend $T$
to a completely positive map on $A + A^*$.
 By Arveson's extension theorem \cite{SOC,Pn}, we may extend further
 to a completely positive map $\tilde{T} : B \to B(H)$.

If  $A$ is an approximately unital
operator algebra, then by Lemmas \ref{refe} and \ref{ik2}, the canonical weak*
continuous extension of $T$ to a map from the unital
operator algebra $A^{**}$ into $B(H)$, is OCP.  By the last paragraph,
 the latter map
has a completely positive
extension $S : B^{**}  \to B(H)$, and $S = V^* \pi(\cdot) V$ for a
$*$-representation $\pi : B^{**}  \to B(K)$ as above.
Restricting $S$ and $\pi$ to $B$ we obtain the desired
extension $\tilde{T} = V^* \pi_{\vert B}(\cdot) V$.

The last assertion,
 about the norm, follows immediately in the unital space case,
since it is well known for
completely positive maps on $C^*$-algebras, and indeed all of our
extensions preserve norms.  If
$A$ is an algebra with cai $(e_t)$, and $B = C^*(A)$,
 then $T(e_t) \to S(1)$
weak*.   Thus $\Vert S(1) \Vert \leq \sup_t \, \Vert T(e_t) \Vert$
by Alaoglu's theorem.  Consequently, by the
unital space case, $\Vert T \Vert_{cb} \leq \Vert S \Vert_{cb}
= \Vert S(1) \Vert  = \Vert V \Vert^2
\leq  \Vert T \Vert$, and so
$\Vert  T \Vert = \Vert T \Vert_{cb} = \sup_t \, \Vert T(e_t) \Vert$.
  \end{proof}

The
next result gives a positive extension into a general $C^*$-algebra,
under a hypothesis that is often satisfied.

\begin{proposition}   \label{ik4}  If $T : A \to B$ is OCP,
from a  unital operator space $A$ into a $C^*$-algebra $B$, and if
there is a (resp.\ weak* continuous) affine map $L : Q(A) \to
Q(C^*(A))$ taking $0$ to $0$, which is a retract of the restriction
map $Q(C^*(A)) \to Q(A)$. Then
 there exists a positive map $\tilde{T} :
C^*(A) \to B^{**}$ (resp.\ $\tilde{T} : C^*(A) \to B$) extending
$T$.
  \end{proposition} \begin{proof}  As before, we
may assume that $C = 1$ in the definition of OCP. If $\varphi \in
S(B)$ then $\varphi \circ T \in Q(A)$. Hence $T^\sharp : Q(B) \to
Q(A) : \varphi \to \varphi \circ T$ is a weak* continuous affine
map.  Then $L \circ T^\sharp : Q(B) \to Q(C^*(A))$ is a (resp.\
weak* continuous) affine map taking $0$ to $0$. For any $c \in
C^*(A)_{\rm sa}$, the map $\epsilon_c : Q(C^*(A)) \to \Cdb$ of
evaluation at $c$, is a  weak* continuous bounded affine map taking
$0$ to $0$. Hence $\epsilon_c \circ L \circ T^\sharp : Q(B) \to
\Cdb$ equals $\epsilon_b$ for a unique $b \in B^{**}_{\rm sa}$
(resp.\ $B_{\rm sa}$), by \cite[3.10.3]{Ped}.  Define $\tilde{T}(c)
= b$. Then $\tilde{T} : C^*(A)_{\rm sa} \to B^{**}_{\rm sa}$ (resp.\
$B_{\rm sa}$) is real linear. Extend $\tilde{T}$ to $C^*(A)$ by
linearity.   If $c \in C^*(A)_+$ then it is clear that
$\psi(\tilde{T}(c)) \geq 0$ for all $\psi \in S(B)$, so $\tilde{T}$
is positive.
\end{proof}

{\bf Remark.}   In the light of the last result, it is worth pointing
out that there need not exist a weak* continuous retract $L : S(A) \to
S(C^*(A))$.  For example, suppose that such a retract existed
when  $A$ is the sum of the compact operators on $\ell^2$
and the upper triangular operators with constant entries on the leading diagonal
(that is, $t_{ij}=0$ unless $j\ge i$, and $t_{ii}=t_{jj}$ for all $i,j$).
The states $\varphi_n$ on $B(\ell^2)$ picking out the $n$th entry on the leading diagonal,
when restricted to $A$, converge weak* on $A$.
However $(L(\varphi_n |_A))$ has no
weak* limit.  Indeed, the restriction of $\varphi_n$ to the compact
operators is well known to have a unique state extension, so $L(\varphi_n |_A) = \varphi_n$.
If  $t=(t_{ij})$ with $t_{ii}=1$ and $t_{2k-1,2k}=1$ ($k\in\bN$),
 and all other $t_{ij}=0$, then the diagonal entries of $tt^*$ are $2,1,2,1,\ldots$ so
$(\phi_n(tt^*))$ does not converge to any limit.

 \bigskip

{\em Acknowledgments:}  We are deeply grateful to the referee for a
superb and incisive effort,
and for very many suggestions and improvements.  In particular, the proof of
Lemma \ref{refe} was generously supplied by the referee, 
answering a conjecture in an earlier version
of the paper (in that version we had proved Theorem \ref{ikc} 
for unital operator spaces,
 and for operator algebras $A$
for which the conclusion of Lemma \ref{refe} was true for all $M_n(A)$).
Lemma \ref{ikhuh} is also due to the referee.  
In several proofs, particularly \ref{abr}, \ref{corea11}, \ref{ruid}, and
\ref{ik3},  we
have followed the referees suggestions for more direct arguments than our original ones.  We also thank Damon Hay for several errata; David Sherman for a discussion that led to us understanding
 the  $C^*$-algebra case of Lemma \ref{refe};  and
 Jean Esterle for independently suggesting to the  first author the example in Section 5,
as a possible candidate for a radical,
approximately unital semiprime operator algebra.

\end{document}

\end{thebibliography}
\end{document}